\documentclass[a4paper]{article}
\usepackage[a4paper,margin=3cm]{geometry}
\usepackage[utf8]{inputenc}
\usepackage{amsthm,amssymb,amsbsy,amsmath,amsfonts,amssymb,amscd,dsfont}
\usepackage{stmaryrd}
\usepackage{xcolor}
\usepackage{graphicx}
\usepackage{placeins}
\usepackage{subcaption}
\usepackage{authblk}
\usepackage{booktabs}
\usepackage{subcaption}
\usepackage{float}
\usepackage{subfiles}
\usepackage{csquotes} 
\usepackage{color}
\usepackage{microtype}
\usepackage{hyperref}
\usepackage{float}
\usepackage{placeins}
\usepackage{float}
\usepackage{subcaption}
\usepackage{booktabs}
\usepackage{tikz}
\usepackage{tikz-cd}
\newtheorem{theorem}{Theorem}[section]
\newtheorem{lemma}[theorem]{Lemma}
\newtheorem{proposition}[theorem]{Proposition}
\newtheorem{corollary}[theorem]{Corollary}

\theoremstyle{definition}
\newtheorem{definition}[theorem]{Definition}

\theoremstyle{remark}
\newtheorem{remark}[theorem]{Remark}
\newcommand{\eps}{\varepsilon}
\newcommand{\R}{\mathbb{R}} 
\newcommand\E{\mathbb E}
\renewcommand\r{\mathbb{R}}

\newcommand\pro[1]{\mathbb{P}\left(#1\right)}
\newcommand\esp[1]{\mathbb{E}\left[#1\right]}
\newcommand\uno[1]{\mathds{1}_{\left\{#1\right\}}}

\title{
Estimation for High-Dimensional Hawkes Processes with Inhibition via CIR Approximation
}

\author[1]{Nadia Belmabrouk}
\author[2]{Mohamed Ben Alaya}
\author[3]{Xavier Erny}

\affil[1]{MODAL'X, Universit\'e Paris Nanterre, Nanterre, 92000, France}
\affil[2]{LMRS, université de Rouen
Avenue de l'Université,
76801 Saint-Étienne-du-Rouvray, France}
\affil[3]{SAMOVAR, T\'el\'ecom SudParis, Institut Polytechnique de Paris, Palaiseau, France}

\date{}

\begin{document}
\maketitle 

\vspace{0.25cm}
\begin{abstract}
In this paper, we develop a parametric estimation method for a high-dimensional Hawkes processes with inhibition. The proposed approach relies on a quantitative approximation of the Hawkes process by a Cox–Ingersoll–Ross (CIR) process. In particular, we establish an explicit convergence rate in the joint large-population and long-time regime. This result requires specific techniques to handle the non-Lipschitz diffusion coefficient of the limiting CIR process.
Based on this approximation, we construct explicit and easily implementable estimators of the model parameters. We establish their consistency and asymptotic normality. We also conduct  numerical experiments to illustrate our findings.
\end{abstract}

\vspace*{0.5cm}


\section{Introduction}

Since their introduction by \cite{hawkes}, Hawkes processes have been widely used to model event dynamics across a broad range of fields. At their core, they provide a framework for describing how the occurrence of an event can alter the rate at which subsequent events occur, thereby capturing interactions between events over time. In the classical formulation, these interactions are excitatory, giving rise to what are commonly referred to as self-exciting processes. Extensions of the framework can also account for inhibitory effects, whereby events decrease the rate of subsequent events. This flexibility has made Hawkes processes useful for modeling a wide variety of phenomena, including earthquake occurrences \cite{hawkes,helmstetter,kagan,ogata}, the co-occurrence of genetic motifs \cite{reynaudbouret}, and interactions in social networks \cite{zhou}. In the following, we focus on two prominent applications: financial order books \cite{AbergelJedidi2015,Hewlett2006,Large2007,LuAbergel2018,MorariuPatrichiPakkanen2022,ProtterWuYang2024} and the activity of biological neural networks \cite{Krumin2010,Lambert2018,Masud2011,ReynaudBouret2014}. These two applications are of particular interest for our model because they involve a large number of interacting components, each generating events at a high frequency. They therefore provide natural examples of systems in which a large number of interactions occur within short time intervals. Since our approach does not require keeping track of individual events, it is well-suited to this context.

In its most general form, a multivariate Hawkes process is a $d$-dimensional point process $(Z^1,...,Z^d)$ with stochastic intensities satisfying: for $1\leq i\leq d,$
\begin{equation}\label{eq:intHawkes}
\lambda^i_t = f_i\left(\sum_{j=1}^d \int_{[0,t)}h_{ji}(t-s) dZ^j_s\right),
\end{equation}
where $f_i$ ($1\leq i\leq d$) are the jump rate functions, and $h_{ji}$ ($1\leq i,j\leq d$) are the interacting functions. The processes $X^i$ defined as
$$X^i_t = \sum_{j=1}^d \int_{[0,t]}h_{ji}(t-s) dZ^j_s$$
are often called potentials. This last term comes from the terminology of neuroscience. In a neural system model, the integer~$d$ represents the number of interacting neurons, each point process~$Z^i$ models the spike trains of each neuron, the process~$X^i$ is the membrane potential of the $i$-th neuron, the function $f_i$ is the firing rate of the neuron depending on its potential, and the functions $h_{ji}$ describe the interaction between the neurons. Formally, the dynamics of the potential~$X^i$ of~\eqref{eq:intHawkes} can be described in two mechanisms:
\begin{itemize}
    \item for each jump time of any~$Z^j$ ($1\leq j\leq d$), every process~$X^i$ ($1\leq i\leq d$) does a jump of height~$h_{ji}(0)$,
    \item between two jumps, the processes~$X^i$ evolve deterministically depending on each past jump event (if $s$ is a previous jump time, it contributes to the dynamics of~$X^i$ at $t\geq s$ through $h_{ji}(t-s)$).
\end{itemize}

In the applications, the jump rate functions~$f_i$ are often assumed to be non-decreasing, such that the variation in $X^i_t$ can be interpreted either as excitation (when $X^i$ increases) or inhibition (when $X^i$ decreases). In particular, if $h_{ji}\geq 0$ (resp. $h_{ji}\leq 0$), the nature of the interaction from~$j$ to~$i$ is excitatory (resp. inhibitory).

The interpretation of~\eqref{eq:intHawkes} in a model of financial order book is similar. The integer~$d$ represents the number of financial actors that can buy or sell some asset, $X^i_t$ corresponding to the tendency for the $i$-th actor to actually buy or sell the asset at the current time~$t$. Naturally, the potentials~$X^i$ depend on the price of the asset, and each time the asset is bought or sold, the price respectively increase or decreases. Hence, if some actor~$j$ is a buyer (resp. seller), the interacting functions of the form~$h_{ji}$ ($1\leq i\leq d$) should be non-negative (resp. non-positive).

In this paper, we work on a mean-field particle system composed of $2N$ particles, divided into two subsets: the exciting particles and the inhibiting particles. In this mean-field framework, the particles share the same potential denoted by~$X^N_t$ with the same rate function~$f$ independent of~$i$, and the exciting (resp. inhibiting) particles~$j$ share the same interacting function $h_{ji} = N^{-1/2} h \geq 0$ (resp. $h_{ji} = -N^{-1/2} h$) for any~$i$. Meaning
$$X^N_t = \frac1{\sqrt{N}}\sum_{j=1}^N \int_{[0,t]} h(t-s)dZ^{N,+,j}_s - \frac1{\sqrt{N}}\sum_{j=1}^N \int_{[0,t]} h(t-s)dZ^{N,-,j}_s,$$
with $Z^{N,+,j},Z^{N,-,j}$ ($1\leq j\leq N$) being orthogonal point processes with the same intensity~$(f(X^N_{t-}))_{t\geq 0}.$ In addition, we assume the convolution kernel to be exponential: for $t\geq 0,$
$$h(t) = e^{-\alpha\, t},$$
where $\alpha>0$ is a model parameter. A classical choice for the rate function~$f$ for Hawkes processes is an affine function (such a model is called a "linear Hawkes process"). However, when we consider Hawkes process with inhibition, it is necessary to enforce the rate function~$f$ to be non-negative. Here we use
$$f(x)= a\, x_+ := a\, \max(x,0),$$
with $a>0$ another model parameter.

Formally, the model we study in this paper is given by
\begin{equation}\label{eq:XN1}
X^N_t = v_0 + (X^N_0-v_0) e^{-\alpha\,t} + \frac1{\sqrt{N}}\sum_{j=1}^N \int_{[0,t]} e^{-\alpha (t-s)}dZ^{N,+,j}_s - \frac1{\sqrt{N}}\sum_{j=1}^N \int_{[0,t]} e^{-\alpha (t-s)}dZ^{N,-,j}_s,
\end{equation}
where the dynamics is characterized by the three positive parameters~$v_0,a,\alpha$. Our main goal concerns the estimation of these parameters by establishing the consistency (see Theorem~\ref{T2}) and the asymptotic normality (see Theorem~\ref{T3}) of the estimators that we introduce in Section~\ref{sec:stat}. The equation~\eqref{eq:XN1} can be classically rewritten as
\begin{align}\label{eq:XN}
    dX^N_t =& \alpha\left(v_0 - X^N_t\right)dt\\
    &+\frac{1}{\sqrt{N}}\sum_{j=1}^N\int_{\mathbb{R}_+}\mathbf{1}_{\left\{z\leq a \left(X^{N}_{t-}\right)_+\right\}}\pi_+^j(dt,dz)-\frac{1}{\sqrt{N}}\sum_{j=1}^N\int_{\mathbb{R}_+}\mathbf{1}_{\left\{z\leq a \left(X^{N}_{t-}\right)_+\right\}}\pi_-^j(dt,dz),\nonumber
\end{align}
where $\pi_-^j,\pi_+^j$ ($j\geq 1$) are independent Poisson measures on $\r_+^2$ with Lebesgue intensity. The scaling factor $1/\sqrt{N}$ corresponds to a diffusive regime. Indeed, since the exciting and inhibiting particles have the same intensity, their contributions cancel at the level of the mean, while their fluctuations accumulate and remain of order one as $N$ tends to infinity. This balance between excitation and inhibition is at the origin of the diffusion limit studied in this paper.

Our model possesses two features that are particularly relevant for neuroscience and order book modeling:
\begin{itemize}
    \item Mean reversion to $v_0$: the drift term $\alpha(v_0-x)$ forces the potential $X^N$ back toward the value $v_0$ at exponential speed, while being perturbed by the jumps. In neuroscience, $v_0$ models the resting potential of the neuron. In finance, if this potential $X^N_t$ is interpreted as a price, $v_0$ represents the long-term equilibrium level ensuring market stability.
    \item Balanced interactions: the excitation and inhibition strengths are identical. This can model an excitation-inhibition balanced neural network, or a balanced order flow in a market order book. Note that, in~\eqref{eq:XN}, we choose the same number of exciting and inhibiting particles with the same interaction strength for excitation and inhibition. The model could however be generalized by considering, for example, a larger number of exciting particles, with a weaker excitation strength so as to preserve the balance of the network. Although such a configuration is also biologically relevant (see, e.g., \cite{BarralReyes2016,LefortEtAl2009}), we choose to focus on the simpler framework with equal numbers of exciting and inhibiting particles. 
\end{itemize}

The estimation of the parameters~$v_0,a,\alpha$ is done asymptotically w.r.t.~$N,t$ going jointly to infinity. Our structure is split into two steps. Firstly, controlling the convergence of the process~$X^N$ to the Cox--Ingersoll--Ross (CIR) process~$\bar X$ in distribution
\begin{equation}\label{eq:Xbar}
d\bar X_t = \alpha\left(v_0 - \bar X_t\right)dt + \sqrt{2a\, \bar X_t}dB_t,
\end{equation}
with $B$ a standard one-dimensional Brownian motion. Secondly, exploiting the results from \cite{alaya2012parameter} to estimate the parameters of our model given by  \eqref{eq:XN}. 

The CIR process, introduced by Cox, Ingersoll and Ross
\cite{CoxIngersollRoss1985}, is a one-dimensional diffusion
process widely used in mathematical finance and stochastic
modeling. In particular, it is used to model interest rates
and stochastic volatility in financial markets
\cite{CoxIngersollRoss1985,Heston1993}. It also finds
applications in biological modeling, particularly in population
dynamics \cite{BansayeMeleard2015}. Under the
assumptions \(v_0,a,\alpha>0\) and a nonnegative initial
condition, equation \eqref{eq:Xbar} admits a unique nonnegative
strong solution. Moreover, the boundary behavior at zero is
governed by the Feller condition. More precisely, if
\(\alpha v_0\geq a\), the process remains strictly positive
whenever it starts from a positive initial condition.
If \(0<\alpha v_0<a\), zero is attainable and instantaneously
reflecting
(see, e.g., \cite{GoingJaeschkeYor2003}).
The CIR process exhibits a mean-reverting behavior and admits
a unique invariant probability distribution, which is a Gamma
distribution with shape parameter \(\alpha v_0/a\) and scale
parameter \(a/\alpha\). Moreover, the process is ergodic
(see, e.g., \cite{Overbeck1998}). 

From a statistical perspective, parameter estimation for the
CIR process has received considerable attention in the
literature. Several approaches have been developed for
estimating the drift and diffusion parameters, including
maximum likelihood methods based on continuous or discrete
observations, as well as moment and least-squares estimation
procedures. The ergodicity of the CIR process plays a
fundamental role in the asymptotic analysis of these
estimators, ensuring consistency and, under suitable
regularity conditions, asymptotic normality
(see, e.g., \cite{OverbeckRyden1997,Overbeck1998}).

More specifically, Ben Alaya and Kebaier investigated
parameter estimation for the CIR process in both
ergodic and nonergodic regimes, establishing asymptotic
properties of the proposed estimators
\cite{alaya2012parameter,BenAlayaKebaier2013}.
Furthermore, Ben Alaya, Kebaier, and Tran established local
asymptotic properties of the CIR model under discrete
observation schemes, providing a detailed analysis of the
corresponding statistical experiments \cite{ben2020local}.

The proofs of our results rely mainly on the convergence of the Hawkes process potential $X^N$ to the CIR process $\bar X$, for which we establish an explicit convergence rate depending on both~$N$ and~$t$. Such quantitative estimates are crucial for our approach, since we consider a joint asymptotic regime in which both~$N$ and~$t$ tend to infinity. Obtaining these estimates raises an important technical difficulty: classical quantitative results for mean-field systems typically rely on (locally) Lipschitz regularity assumptions on the coefficients (see, e.g., \cite{JourdainMeleard1998,Sznitman1991}). In our setting, however, the diffusion coefficient in \eqref{eq:Xbar}, which involves the square-root function, is not Lipschitz continuous at zero. We overcome this difficulty by introducing a smooth cut-off of the square-root function around zero (see Section~\ref{sec:31}). While standard Lipschitz estimates typically lead to a Gr\"onwall-type factor of order $N^{-1/2}e^{C\,t}$, our approach, in this non-locally Lipschitz framework, yields a convergence rate of order $N^{-1/2}e^{t/\eps_N} + t\eps_N$, where $\eps_N$ is any vanishing positive sequence (see Corollary~\ref{cor:XNbarX} for a formal statement). The joint asymptotic regime considered in our statistical analysis is chosen so that both terms in this bound vanish. In comparison, in \cite{erny2022mean}, the authors study a similar model, where the jump rate function~$f$ is generic, and the diffusive coefficient of the limit SDE $\sqrt{f}$ is assumed to be Lipschitz continuous. In that work, the regularity assumption on the limiting diffusion coefficient allows one to rely on standard Lipschitz estimates.

Scaling limits of Hawkes processes with non-exponential convolution kernels have also been investigated (see in particular \cite{AbiJaberCuchieroLarssonPulido2021}). In this setting, the limiting dynamics can be described by stochastic Volterra equations, for which a weak solution theory is developed in that reference. The corresponding limit results are qualitative and do not provide an explicit convergence rate, which is essential for our analysis.

To study our estimators, we also need to consider some functionals of $X^N$, in addition to the process itself. While $X^N$ is a time-homogeneous Markov process, the process obtained by adding these functionals to $X^N$ is a time-inhomogeneous Markov process. We can therefore use its generator to establish the required convergence results.

Regarding parameter estimation, maximum likelihood estimation (MLE) is one of the most common approaches in the Hawkes process literature. Due to the complexity of the likelihood function, numerical methods are generally used to compute the MLE; see, for instance, \cite{bonnet2023inference} and \cite{LuAbergel2018} which apply numerical MLE methods to high frequency Hawkes models with excitation and inhibition in order book data.

In practice, our estimators depend exclusively on the values of the process~$(X^N_t)_{t\geq 0}$, without requiring explicit knowledge of its jump times, unlike the maximum likelihood estimator in this framework. Although the jump times can theoretically be recovered from the continuous observation of $X^N$, this information may not be readily available in applications. Indeed, in high-frequency applications (i.e., when $N$ is large), the process $X^N$ behaves similarly to the continuous process~$\bar X$, which exhibits a diffusive behavior (see Figure~\ref{figg}). Our estimators, therefore, have the advantage of relying only on the observed values of the process, without requiring the identification of individual jumps. Moreover, our estimators are computationally simpler, and their computational complexity does not depend on the number of particles, while avoiding the convexity issues associated with maximum likelihood estimation.

The paper is structured as follows: Section \ref{sec:stat} introduces the statistical problem and the associated estimation problem. Section \ref{sec:main} presents the main results. Section \ref{sec:num} is devoted to numerical experiments. The proofs of Theorem~\ref{T1}, \ref{T2}, and \ref{T3} are provided in Sections \ref{sec:31}, \ref{sec:33}, and \ref{sec:34}, respectively. Finally, Sections \ref{sec:A}, \ref{sec:B}, and \ref{sec:C} contain additional technical details and supplementary results.

\section{Statistical model}\label{sec:stat}
We assume that we observe continuously the potential of the Hawkes process~$(X^N_t)_{0\leq t\leq T}$ defined at~\eqref{eq:XN}, where the unknown parameters are~$v_0,\alpha,a$, that are assumed to be positive. The time-horizon~$T>0$, the number of exciting (resp. inhibiting) particles~$N\in\mathbb{N}^*$, and the initial condition~$X^N_0=x>0$ are known. The estimation is done asymptotically w.r.t.~$(N,T)$ going jointly to infinity.

\begin{proposition}
The statistical model consisting of the observation of~$(X^N_t)_{0\leq t\leq T}$ and of the unkown parameter~$(v_0,\alpha,a)$ is identifiable.
\end{proposition}

\begin{proof}
Firstly, the identification of~$a$ comes from the fact that
$$\esp{[X^N]_t} = 2a\int_0^t \esp{(X^N_s)_+}ds,$$
hence, denoting $x=X^N_0>0,$
$$a = \frac{1}{2x}\frac{d}{dt}\left(\esp{[X^N]_t}\right)_{t=0}.$$

In order to identify~$v_0$ and~$\alpha$, let us notice that, by denoting $A^N$ the generator of~$X^N$ (in the sense of Definition~\ref{def:semi} for example), for any reachable~$y>0$ by the process~$X^N,$
$$A^N(\textrm{Id})(y) = \alpha(v_0-y).$$

Whence, it is possible to indentify $\alpha$ and $v_0$ by choosing two distinct values~$y_1,y_2$ different from~$v_0$.
\end{proof}

In order to study the parameters, we consider the two following processes, defined as functionals of $X^N$:
$$I^N_t= \frac{1}{t} \int_{0}^{t}  X^N_s \, ds~~\textrm{ and }~~
J^N_t= \frac{1}{t} \int_{0}^{t} \frac{1}{ X^N_s +\eta_N}\, ds,$$
where $\eta_N$ is a positive sequence verifying
$$N^{-1/2}<< \eta_N<<1.$$

Unlike the CIR process~$\bar X$ that cannot reach zero under Feller condition, the process~$X^N$ can take negative values. In fact, it can be proved that, for any~$t\geq 0$, $X^N_t > -N^{-1/2}$. This is the reason why, in the definition of~$J^N$ above, we use the sequence~$\eta_N$.

Based on the observations of the process $X^N$, we define the following parameter estimators: 
\begin{enumerate}
    \item $\displaystyle\widehat{v}_0^{N}(T)=\frac{1}{T}\int_0^T X^N_s\,ds=I_T^N.$
    \item $\displaystyle\widehat{a}^{N}(T)=\frac{1}{2}\frac{[X^N]_T}{\int_0^T X_s^N\,ds}=\frac{[X^N]_T}{2TI_T^N}.$
    \item $\displaystyle\widehat{\alpha}^N(T)=\frac{\hat a^N(T)}{\hat v_0^N(T) - (1/J^N_T)},$
\end{enumerate}
where we recall that
\begin{equation}\label{eq:crochet}
[X^N]_t = \underset{p\to 0}{\lim} \sum_{k=0}^{n-1} \left(X^N_{t^p_{k+1}} - X^N_{t^p_k}\right)^2
\end{equation}
in probability, with $0=t^p_0<t^p_1<...<t^p_n=t$ any subdivision of mesh size~$p$.

{Notice that we choose the estimators for~$v_0$ and~$\alpha$ from the results of \cite{alaya2012parameter}, where the authors observe the CIR process~\eqref{eq:Xbar} over a long time horizon. Hence, these estimators are based on the method of moments, whose principle is to match theoretical moments with their empirical counterparts. Here, the theoretical moments are those of the limiting distribution of $\bar X$. Indeed, $\bar X$ is ergodic and its limiting law is the Gamma distribution with shape parameter $\alpha v_0 / a$ and rate parameter $\alpha / a$. The ergodic theorem therefore yields, for any function $f$ integrable under this law,
\[
\frac{1}{T}\int_0^T f(\bar X_s)\,ds \;\xrightarrow[T\to\infty]{\text{a.s.}}\; \E[f(\bar X_\infty)].
\]
Applying this with $f(x) = x$ and, when $\alpha v_0 > a$, with $f(x) = 1/x$, the expectations $\E[\bar X_\infty]$ and $\E[1/\bar X_\infty]$ are consistently approximated by the empirical means
$
\frac{1}{T}\int_0^T \bar X_s\,ds$
and
$
\frac{1}{T}\int_0^T \frac{1}{\bar X_s}\,ds.
$
For this Gamma distribution, these expectations are given by
$
\E[\bar X_\infty] = v_0
$
and
$
\E\left[\frac{1}{\bar X_\infty}\right] = \frac{\alpha}{\alpha v_0 - a},
$
so that the parameters can be expressed as
\[
v_0 = \E[\bar X_\infty]
\qquad \text{and} \qquad
\alpha = \frac{a}{\E[\bar X_\infty] - 1/\E\left[\dfrac{1}{\bar X_\infty}\right]}.
\]
These identities express $v_0$ and $\alpha$ as continuous functions of the two expectations. Replacing the latter by their empirical counterparts therefore gives estimators of $v_0$ and $\alpha$ which are consistent by the continuous mapping theorem.}

In addition, if we observe the process~$\bar X$, then the estimation for~$a$ would be "perfect" in the sense that, for all~$t>0,$
$$\frac{[\bar X]_t}{2 \int_0^t \bar X_s\,ds} = a.$$

Indeed, since~$\bar X$ is a continuous semimartingale, its quadratic variation coincides with its predictable quadratic variation, which is given by
$$\langle \bar X\rangle_t=2a\int_0^t\bar X_sds.$$

For the process~$X^N$, the situation is more intricate since it is not continuous, and hence~$\widehat a^N(T)$ is not equal to~$a$ in general. More precisely, the quadratic variation~$[Y]$ of a semimartingale~$Y$ is determined by its trajectory (see~\eqref{eq:crochet}). 
 In contrast, the predictable quadratic variation~$\langle Y\rangle$ (more precisely, that of its local martingale part) is the predictable compensator of its quadratic variation, and therefore depends on the underlying probabilistic structure of the process. In our statistical framework, the quadratic variation~$[X^N]_t$ is observable, whereas its predictable quadratic variation~$\langle X^N\rangle_t$ is not. The two quantities coincide for continuous semimartingales, such as the CIR process~$\bar X$, which explains why the estimation of~$a$ would be perfect in that case.

We state our statistical results about this model in the next section. Namely, the triple $(\widehat{v}_0^{N}(T),\widehat{a}^{N}(T), \widehat{\alpha}^{N}(T))$ is both consistent and asymptotically normal w.r.t.~$(v_0,a,\alpha)$.

\section{Main results}\label{sec:main}

The three estimators defined in the previous section rely on the processes~$I^N,J^N$. The goal of the next result is to prove that $(I^N_T,J^N_T)$ converges in distribution to $(\bar I_T,\bar J_T)$ as $N$ goes to infinity, with an explicit convergence speed depending on~$(N,T)$. Having an explicit convergence rate is crucial for our approach since we use Theorem~\ref{T1} below by letting $N,T$ go to infinity jointly. Notice that, to prove this result, it is required to consider the triple $(X^N,I^N,J^N)$ to treat a (non-homogeneous) Markov process.

\begin{theorem}\label{T1}
Let $(X_t^N)_{t\geq 0}$ be the solution to~\eqref{eq:XN} and let $(\bar{X}_t)_{t\geq 0}$ be the CIR process solving~\eqref{eq:Xbar}. Assume that $4a<\alpha\,v_0$. Then, for $g:\r^3\to\r$ of class~$C^3_b$, and $x>0,$
$$\left|\mathbb{E}_x\left[g\left(X_T^N,I_T^N,J_T^N\right)\right] - \mathbb{E}_x\left[g\left(\bar X_T,\bar I_T,\bar J_T\right)\right]\right|{\longrightarrow} 0,$$
as $N,T$ tend jointly to infinity, as soon as there exists $M_T\to\infty$ such that $N_T\geq e^{T^2 M_T}.$
\end{theorem}

The previous result gives, in particular, the convergence of the uni-dimensional distributions of the (potential of the) Hawkes process~$X^N$. 
 Using that it is a Markov process, it can be proved by induction that the finite-dimensional distributions of $X^N$ also converge to those of the CIR process~$\bar X$. Together with a tightness argument (whose proof is classical from Aldous' criterion), this leads to the next result. Notice that this kind of quantitative conditional propagation of chaos is not usual for non locally Lipschitz coefficients. Instead of retrieving the classical convergence rate of the form $N^{-1/2}e^{C\,T}$, we obtain one depending on a free parameter~$\eps$ in the statement of Corollary~\ref{cor:XNbarX} below.

\begin{corollary}\label{cor:XNbarX}
Assume that $4a<\alpha\,v_0$. The process~$X^N$ converges in distribution to~$\bar X$ on Skorohod space. In addition, there exists~$C>0,$ such that, for all $g\in C^3_b(\mathbb{R}), N\in\mathbb{N}^*,T>0,x>0,\eps>0$
$$\left|\mathbb{E}_x\left[g\left(X^N_{T}\right)\right] - \mathbb{E}_x\left[g\left(\bar X_{T}\right)\right]\right|\leq C (1+x+1/x)||g||_{3,\infty}\left(N^{-1/2}\, e^{T/\eps} + T\,\eps\right),$$
with
$$||g||_{3,\infty}=\sum_{k=0}^3 ||g^{(k)}||_\infty.$$
\end{corollary}

For example, one can easily use Corollary~\ref{cor:XNbarX} to prove the convergence $X^N$ to $\bar X$ on any compact time interval with an explicit rate, by choosing wisely~$\eps$ as function of~$N$ (e.g. $\eps_N = 1/\ln \ln N$).
\begin{figure}[H]
\centering
\includegraphics[width=0.6\textwidth]{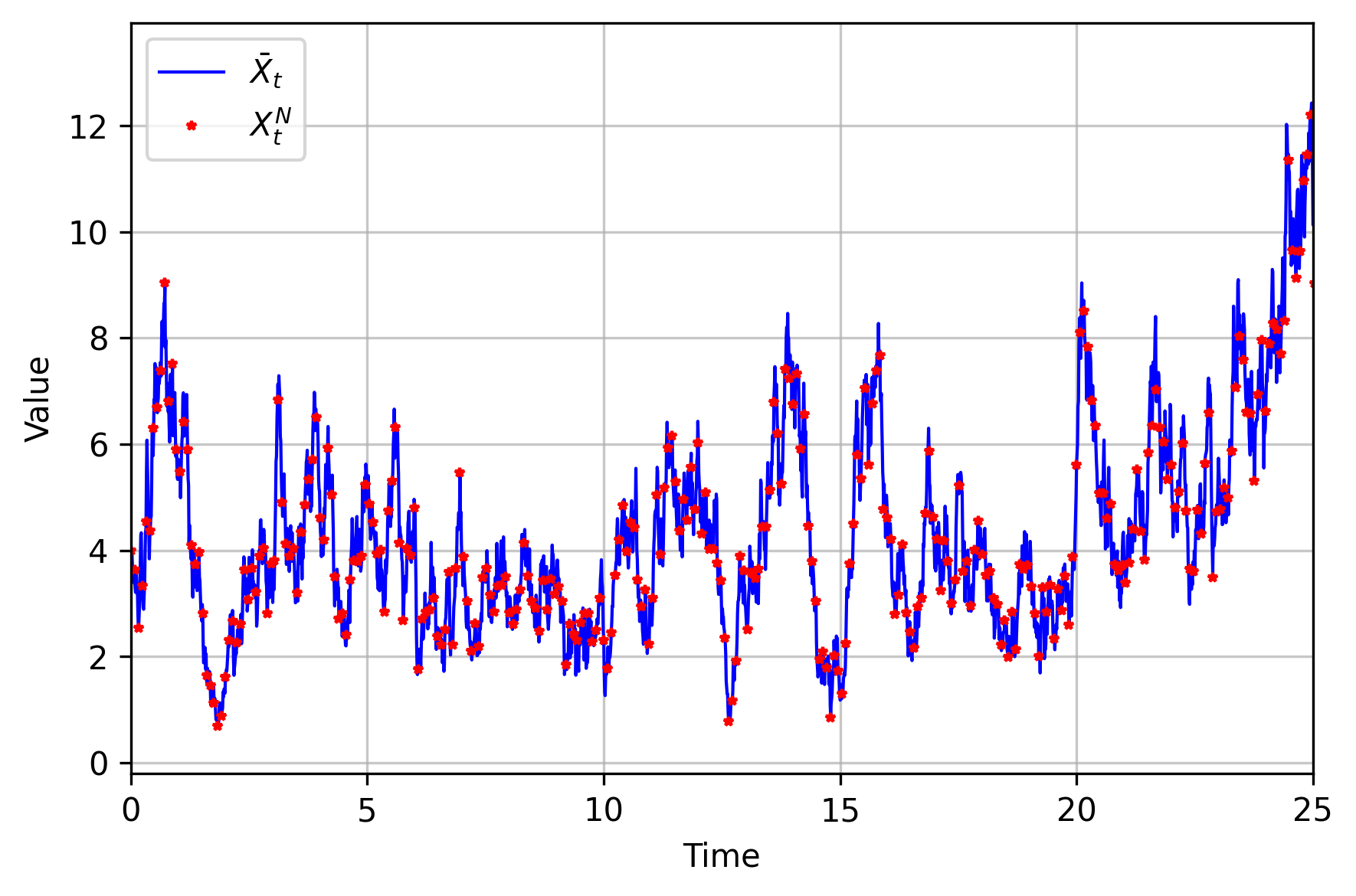}
\caption{Coupling of $X^N_t$ and $\bar{X}_t,$  with $v_0=4,\, \alpha=2,\, a=1.5.$}
\label{figg}
\end{figure}
Figure~\ref{figg} illustrates a numerical approximation of the potential of the Hawkes~$X^N$ process by the CIR process~$\bar X$.  The coupling between the jump process $X^N$ and $\bar{X},$  is constructed following an Euler-Maruyama scheme as follows: at each time step $t_k,$ we generate a common random variable $U_k \sim \mathcal{U}(0,1)$. We then apply the inverse cumulative distribution function of the Skellam distribution to $U_k$ to generate the increment of $X^N,$ and the inverse cumulative distribution function of the standard normal distribution to the same random variable $U_k$ to generate the corresponding increment of $\bar{X}.$

Independently of Corollary~\ref{cor:XNbarX}, we use Theorem~\ref{T1} to study asymptotic properties of the estimators defined in the previous section. The first main result is the consistency of the estimators.

\begin{theorem}[Consistency of the estimators]\label{T2}
Under the condition $4a<\alpha\,v_0$, the estimators
$\widehat{v}_0^{N}(T)$, $\widehat{\alpha}^{N}(T)$, and $\widehat{a}^{N}(T)$
are consistent as $t,N\to\infty$. More precisely, they converge in probability respectively to~$v_0,\alpha,a$ as $N,T$ tend jointly to infinity, as soon as there exists $M_T\to\infty$ such that $N_T\geq e^{T^2 M_T}.$
\end{theorem}


The second main result is the asymptotic normality of these estimators. Note that, for the parameter~$a$, the limit normal distribution is degenerated with a null variance.

\begin{theorem}\label{T3}
Let $(G_I,G_J)$ be a centered Gaussian vector with the following covariance matrix
$$\begin{pmatrix}
\displaystyle\frac{2av_0}{\alpha^2}
&
\displaystyle-\frac{2a}{\alpha(\alpha v_0-a)}
\\[3mm]
\displaystyle-\frac{2a}{\alpha(\alpha v_0-a)}
&
\displaystyle\frac{2a\alpha}{(\alpha v_0-a)^3}
\end{pmatrix}.$$

Under the condition $4a<\alpha\,v_0$, the triple of estimators
$(\widehat{v}_0^{N}(T),\widehat{a}^{N}(T), \widehat{\alpha}^{N}(T))$
satisfies an asymptotic normality property as $T,N\to\infty$. More precisely,
$$\left(\sqrt{T}\left(\widehat{v}_0^{N}(T)-v_0\right),\sqrt{T}\left(\widehat{a}^{N}(T)-a\right), \sqrt{T}\left(\widehat{\alpha}^{N}(T)-\alpha\right)\right)\longrightarrow \begin{pmatrix}
1&0&0
&
\\[3mm]
0&0&0
\\[3mm]
\displaystyle-\frac{\alpha^2}{a}
&\displaystyle \frac{\alpha}{a}
&\displaystyle -\frac{a}{(v_0-a/\alpha)^2}
\end{pmatrix}\begin{pmatrix}
G_I\\[3mm]
0\\[3mm]
G_J
\end{pmatrix}$$
in law as $N,T$ tend jointly to infinity, as soon as there exists $M_T\to\infty$ such that $N_T\geq e^{T^{7/2} M_T}.$
\end{theorem}

The degeneracy w.r.t.~$a$ is due to the fact the fluctuation of~$\widehat a^N(T)$ around~$a$ are not of order~$\sqrt{T}$, but rather~$\sqrt{N\,T}$. While our approach does not allow to prove the asymptotic normality in this regime, the following technical remark can easily be deduced from it.
\begin{remark}\label{rem:flucta}
    There exists a family of r.v. $W(N,T)$ such that
    \begin{equation}\label{eq:t3a2}
    \esp{W(N,T)^2}\underset{N,T\to\infty}{\longrightarrow} 2a\,v_0
    \end{equation}
    such that
    $$\sqrt{N\,T}\left(\widehat a^N(T) - a\right) - \frac1{2v_0}W(N,T)\xrightarrow[T,N\to\infty]{\mathbb{P}} 0$$
    in the same convergence regime as in Theorem~\ref{T3}. This result comes from the calculation of Section~\ref{sec:consa}. Indeed, with the notation introduced therein,
    $$W(N,T)= \sqrt{N} \frac{M^N_T}{\sqrt{T}},$$
    and
    \begin{equation}\label{eq:t3a3}
    \sqrt{N\,T}\left(\widehat a^N(T) - a\right) - \frac1{2v_0}W(N,T) = W(N,T)\frac12\left(\frac{1}{I^N_T} - \frac1{v_0}\right) + a\sqrt{N\,T}\left(\frac{\tilde I^N_T}{I^N_T}-1\right).
    \end{equation}
    Firstly, the convergence of the second-order moment of~$W(N,T)$ (namely~\eqref{eq:t3a2}) is obtained by~\eqref{a1}. This implies that the family $(W(N,T))_{N,T}$ is bounded in probability. Hence, the first term of the RHS of~\eqref{eq:t3a3} vanishes in probability by~\eqref{a0}. The same holds for the second one, as a consequence of the uniformity with respect to~$\eta$ of the bound~\eqref{eq:erra2}. 
\end{remark}

Notice that, as an elementary result, it is possible to calculate directly the bias of $\widehat{v}_0^N(T)$.  
Indeed, by Lemma~\ref{L1}$.(ii)$,
$$\esp{I^N_T} = \frac1T\int_0^T \esp{X^N_s}ds = \frac1T\int_0^T \left[v_0 + (x-v_0)e^{-\alpha s}\right]ds = v_0 + \frac{x-v_0}{\alpha\,T}\left(1-e^{-\alpha T}\right),$$
implying that the bias is vanishing as $T$ goes to infinity uniformly w.r.t.~$N$:
$$\esp{\widehat{v}^N_0(T)} - v_0 = \frac{x-v_0}{\alpha\,T}\left(1-e^{-\alpha T}\right).$$

While the theoretical results stated in Theorems~\ref{T1},~\ref{T2} and~\ref{T3} rely on a restrictive condition on the joint asymptotic regime~$(N,T)\to(\infty,\infty)$, the numerical simulations presented in the next section are not subject to such a restriction.

\section{Numerical experiments}\label{sec:num}
For a given phenomenon, for instance in neuroscience or in the context of limit order books, the process $X^N$ defined in \eqref{eq:XN} can be observed at high frequency in two different ways in practice: either through its trajectory or through its jump times and jump types.
We organize the numerical study into three parts. First, we investigate the finite-sample performance of the proposed estimators based on high-frequency trajectory observations. Second, we present the numerical performance of the maximum likelihood estimator (MLE) based on jump-time and jump-type observations, under the same simulation framework. Finally, we discuss the two estimation procedures associated with the two observation schemes.
We assess the finite-sample performance of the two estimators with a Monte-Carlo
study. The process $X^N,$ which solves equation \eqref{eq:XN} is simulated as follows: the jump times are first generated using Ogata's thinning method \cite{ogata1981lewis} and between two consecutive jump times $t_k$ and $t_{k+1},$ the process evolves according to the explicit expression $X_t^N = v_0 + (X_{t_k}^N-v_0)e^{-\alpha(t-t_k)}.$ We consider the parameter values $v_0=4,\, \alpha=2$ and $a=1.5,$ satisfying the condition $4a<\alpha v_0$. The initial value is set to $X_0 = 5$, which is intentionally distinct from the true value $v_0.$ The
numerical procedure is implemented in Python. 
\subsection{Proposed estimators}
\begin{figure}[H]
\centering
\begin{minipage}{0.32\textwidth}
    \centering
    \textbf{(i)}\\
    \includegraphics[width=\linewidth]{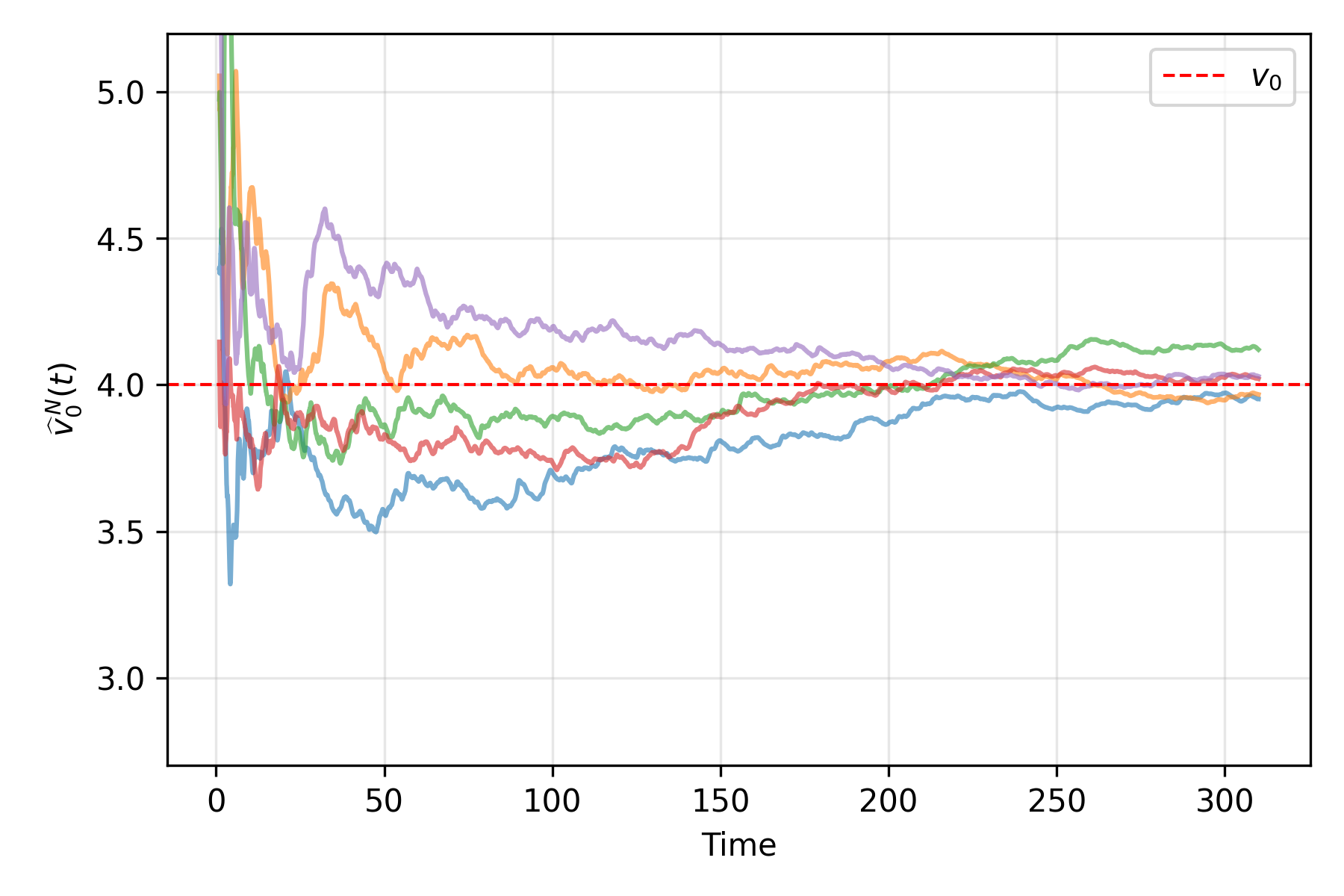}
\end{minipage}
\hfill
\begin{minipage}{0.32\textwidth}
    \centering
    \textbf{(ii)}\\
    \includegraphics[width=\linewidth]{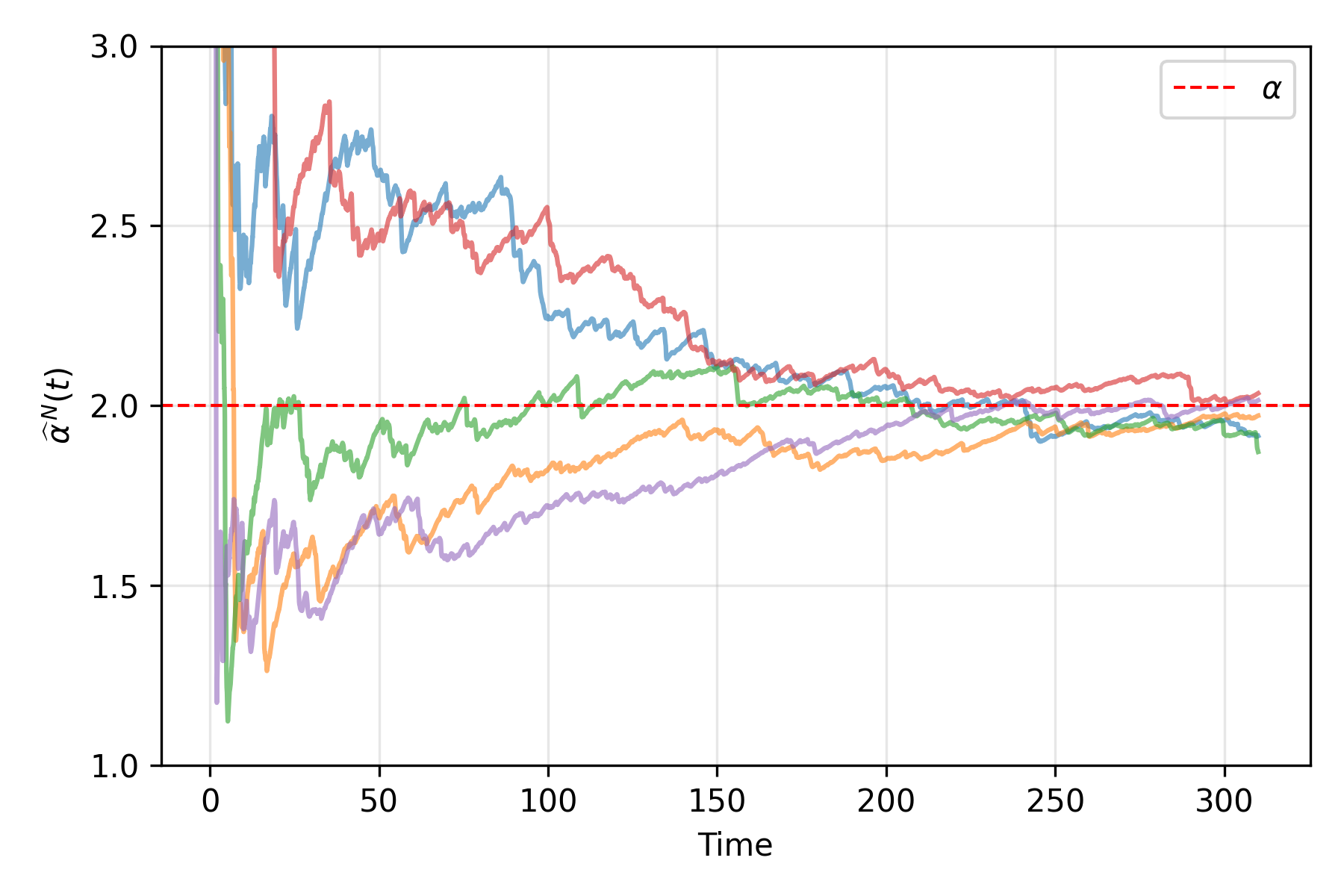}
\end{minipage}
\hfill
\begin{minipage}{0.32\textwidth}
    \centering
    \textbf{(iii)}\\
    \includegraphics[width=\linewidth]{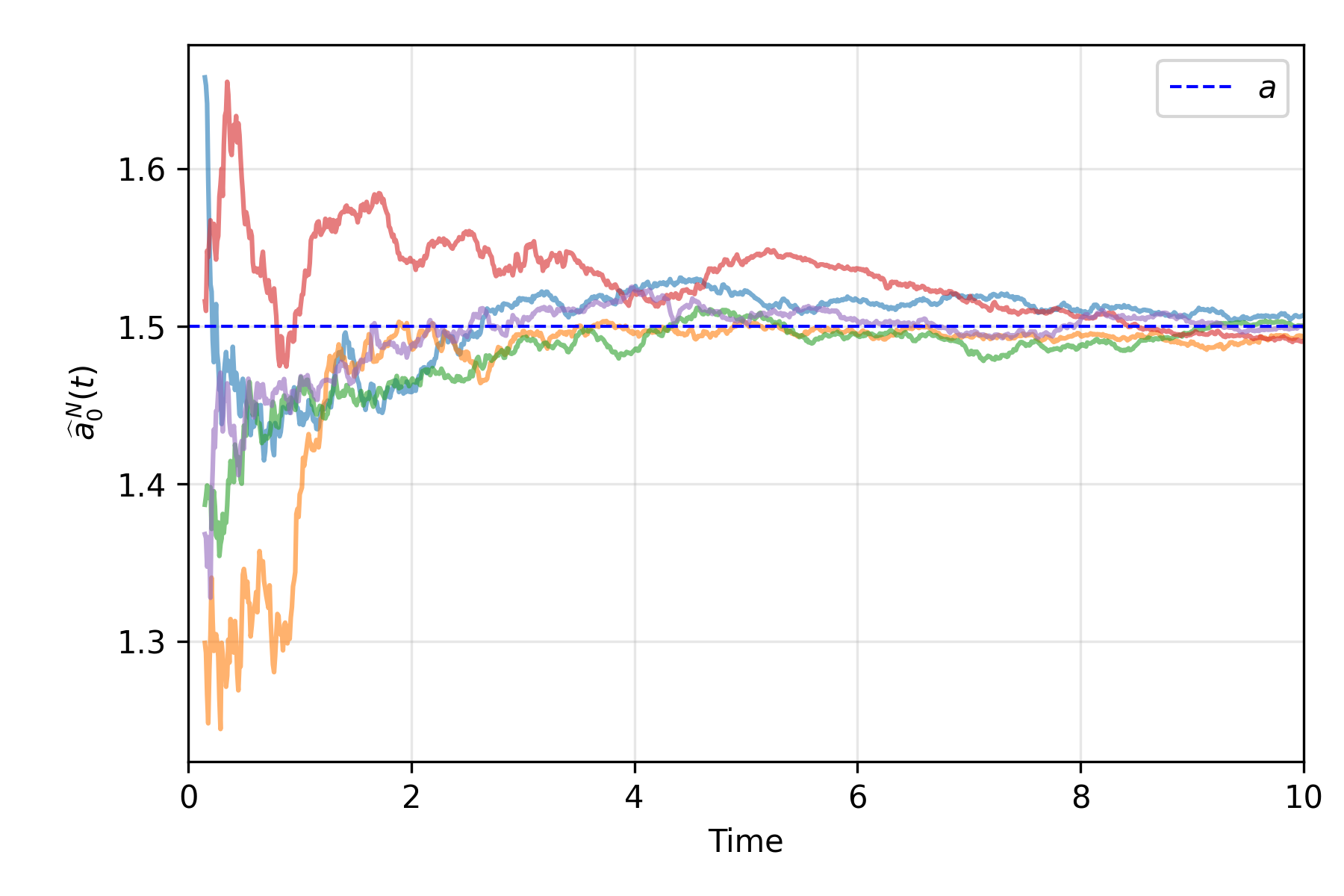}
\end{minipage}

\caption{Consistency of the proposed estimators of
$v_0$ (i), $\alpha$ (ii), and $a$ (iii) over 5 simulated trajectories, with
$v_0=4$, $\alpha=2$, $a=1.5$, and $N=1500$.}
\label{fig:consistency}
\end{figure}

\begin{figure}[H]
\centering

\begin{minipage}{0.32\textwidth}
    \centering
    \includegraphics[width=\linewidth]{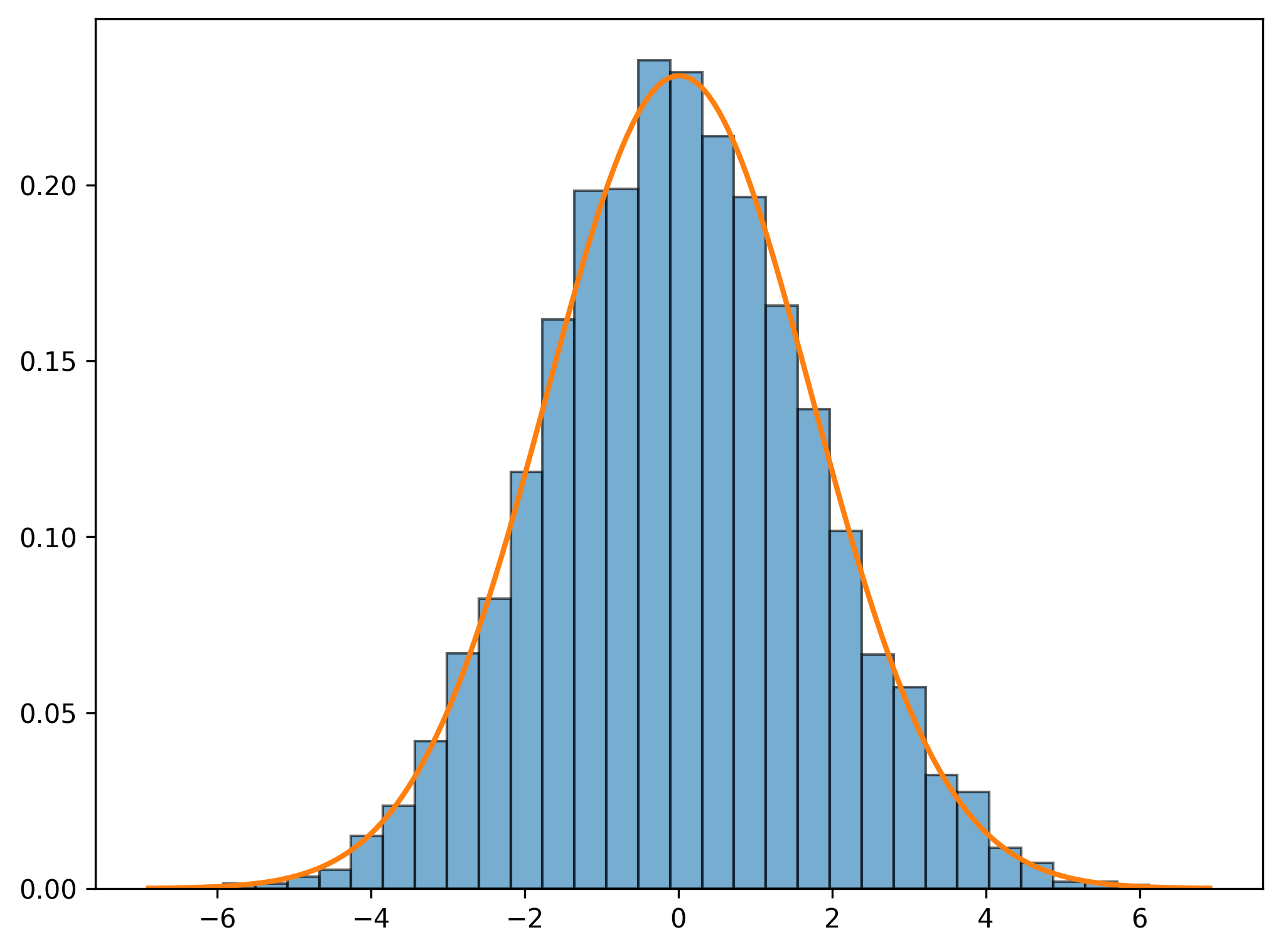}

    \small $\sqrt{T}\left(\widehat{v}_0^N(T)-v_0\right)$
\end{minipage}
\hfill
\begin{minipage}{0.32\textwidth}
    \centering
    \includegraphics[width=\linewidth]{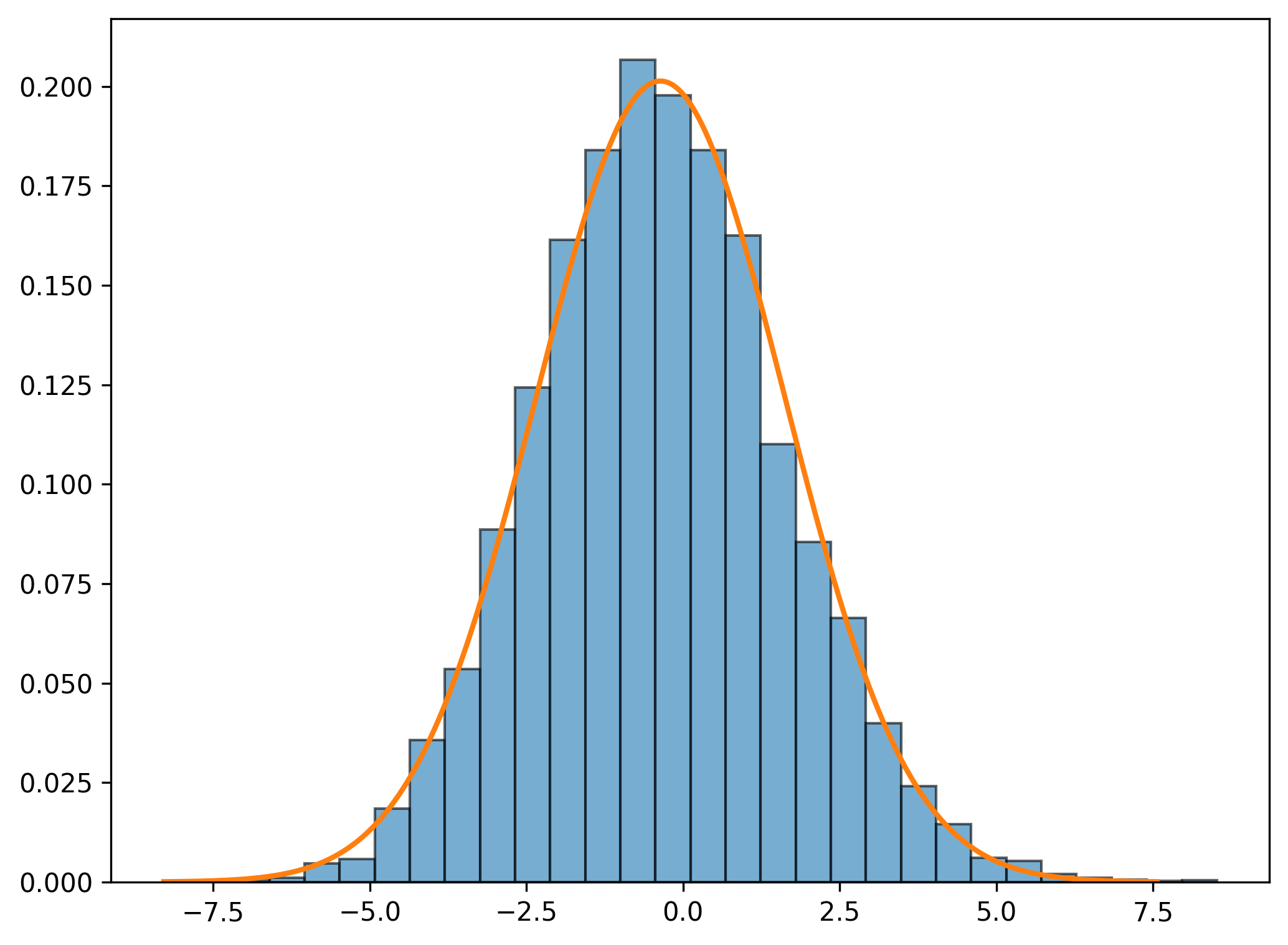}

    \small $\sqrt{T}\left(\widehat{\alpha}^N(T)-\alpha\right)$
\end{minipage}
\hfill
\begin{minipage}{0.32\textwidth}
    \centering
    \includegraphics[width=\linewidth]{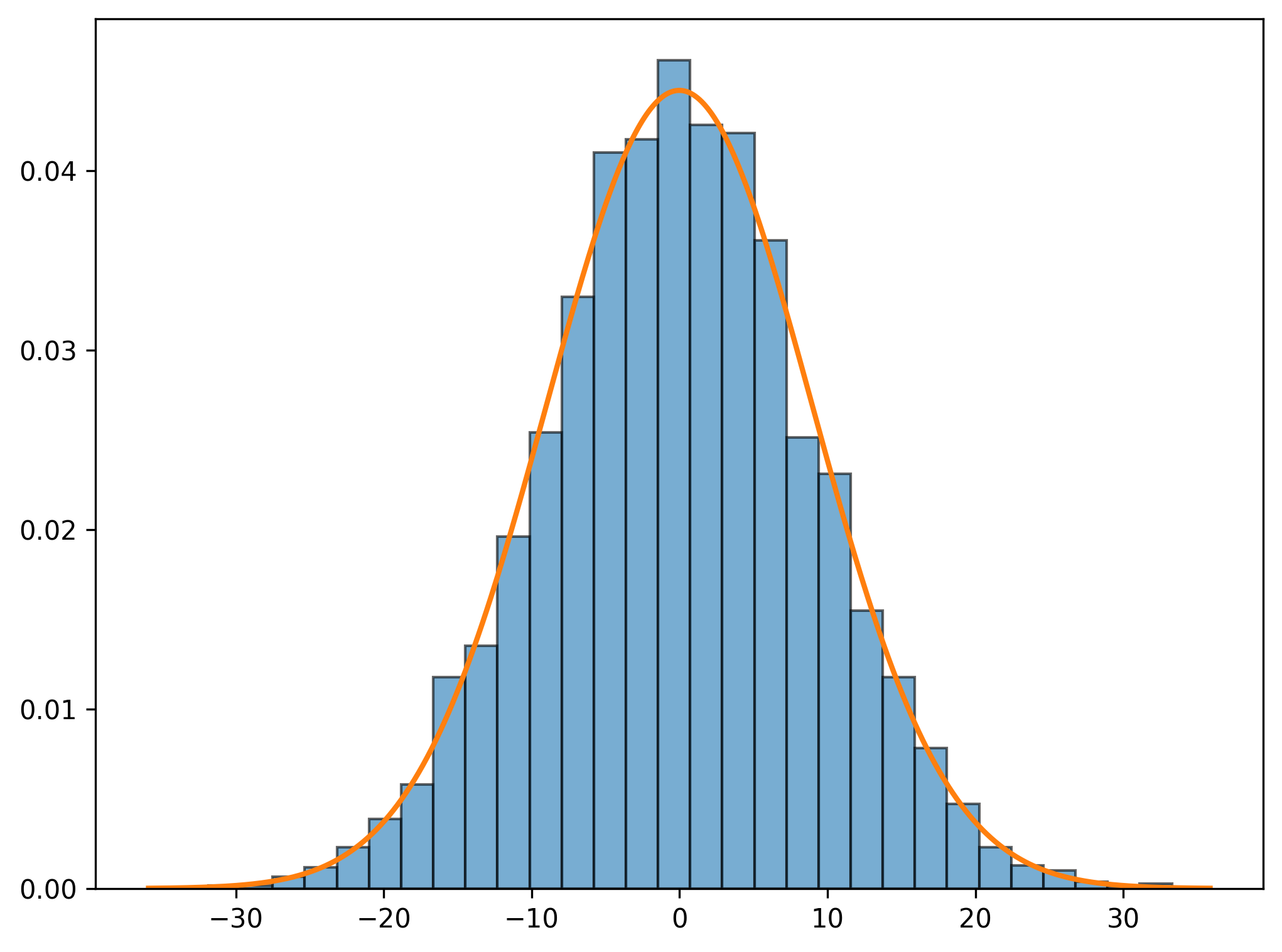}

    \small $\sqrt{NT}\left(\widehat{a}^N(T)-a\right)$
\end{minipage}

\caption{Gaussian fluctuations of the proposed estimators of
$v_0$, $\alpha$, and $a$ around their true values, with
$v_0=4$, $\alpha=2$, $a=1.5$, $T=300$, and $N=1500$,
$M=8000$ Monte Carlo replications.}

\label{fig:fluctuations}
\end{figure}

\begin{figure}[H]
    \centering
    \includegraphics[width=0.8\textwidth]{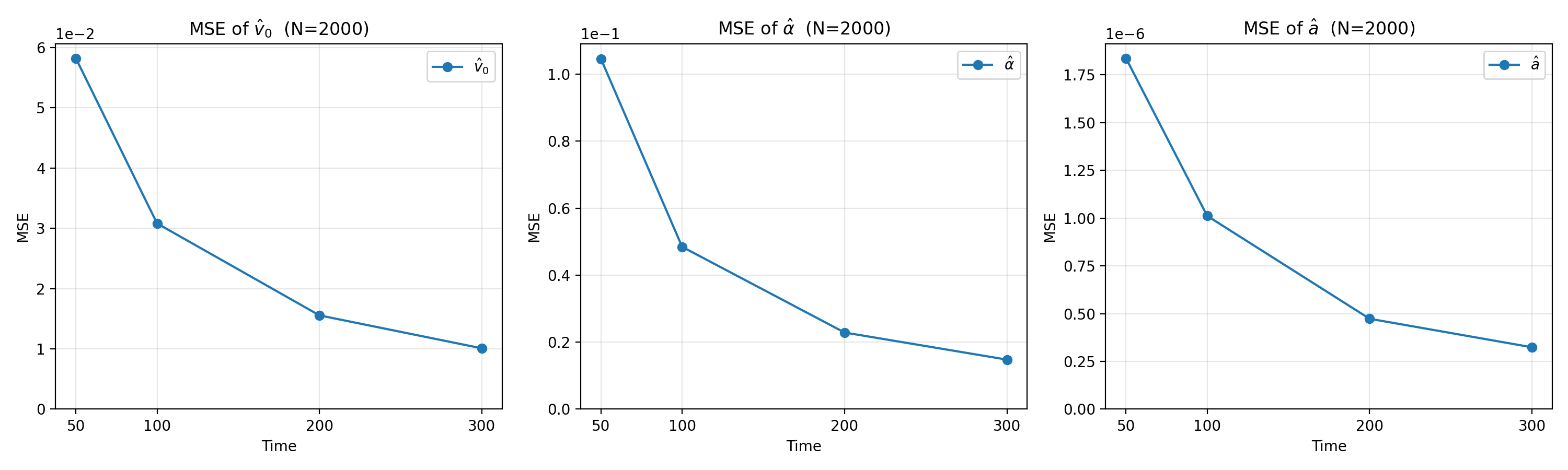}
    \caption{MSEs of the estimated parameters over time with
$v_0=4$, $\alpha=2$, $a=1.5$, $N=2000$ and $M=500$ Monte Carlo replications.}
    \label{MSE4:graph}
\end{figure}

\begin{figure}[H]
    \centering
    \includegraphics[width=0.8\textwidth]{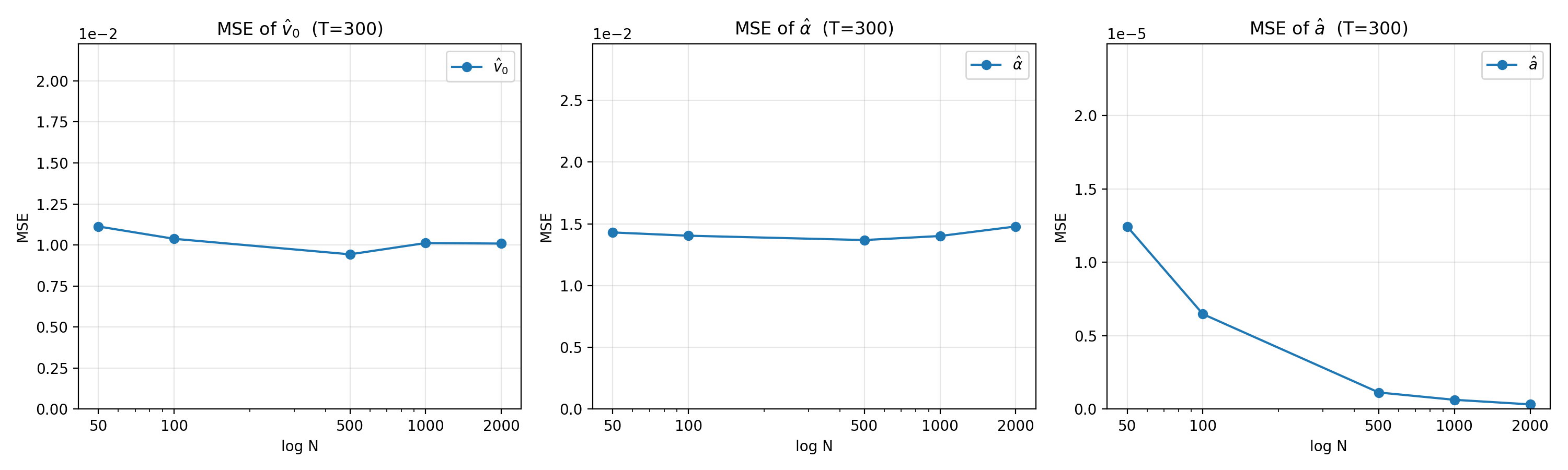}
    \caption{MSE of the Estimated Parameters over $N$ with
$v_0=4$, $\alpha=2$, $a=1.5$, $T=300$ and $M=500$ Monte Carlo replications.}
    \label{MSE5:graph}
\end{figure}

 The estimators are obtained analytically in closed form using a discretized trajectory. The integrals $\int_0^T X^N_sds $ and $\int_0^T \frac{1}{X^N_s + \eta_N}ds$ are approximated by Riemann sums with time step $\Delta_n = 0.01,$ and the quadratic variation $[X^N]_T$ approximated by the sum of the squared increments of the discretized path.
Figure \ref{fig:consistency} illustrates the consistency of the proposed estimators. Each panel displays five simulated trajectories of the corresponding estimator, which converge toward the true parameter value. Figure \ref{fig:fluctuations} displays histograms of the normalized estimation errors obtained from $8000$ Monte Carlo replications.
We observe that the empirical distributions are approximately centered at zero and close to a Gaussian shape. Figure \ref{MSE4:graph} shows that the MSEs of the estimated parameters decrease as $T$ increases, indicating that estimation accuracy improves with longer observation periods. Figure \ref{MSE5:graph} shows that the MSE of the estimator of $a$ decreases as the sample size $N$ increases, whereas the MSEs of the estimators of $v_0$ and $\alpha$ remain almost unchanged. These numerical results confirm  the theoretical findings established in Theorem~\ref{T3} and Remark~\ref{rem:flucta}: the fluctuations related to $\hat v^N_0(T)$ and $\hat \alpha^N(T)$ (resp. $\hat a^N(T)$) are of order~$T^{-1/2}$ (resp. $(NT)^{-1/2}$).

\subsection{Comparison with Maximum likelihood estimator }
Consider $\Theta=\{(v_0,\alpha,a)\in (0,+\infty)^3:\, 4a<\alpha\, v_0\}$. Suppose now that we observe the event times associated with the excitation and inhibition processes over the interval $(0,T],$ namely $\{T_k^{N,+} : 1 \leq k \leq Z^{N,+}(T)\} 
\quad \text{and} \quad
\{T_k^{N,-} : 1 \leq k \leq Z^{N,-}(T)\},$ where $Z^{N,+}$ and $Z^{N,-}$ denote the counting processes associated with excitation and inhibition, respectively. 
For each $\theta \in \Theta$, the log-likelihood function is given by (see Proposition 7.3(III) in \cite{daley2003introduction})
\begin{equation}\label{MLE}
    \ell_T(\theta)=\ell_T ^+(\theta)+\ell_T^-(\theta),
\end{equation}
where 
$$ \ell_T^+(\theta) = \sum_{k=1}^{Z^{N,+}(T)} \log\lambda_\theta^+(T_k^{N,+}) -\Lambda_\theta^+(T),$$
and 
$$\ell_T^-(\theta)=\sum_{k=1}^{Z^{N,-}(T)}\log\lambda^-_{\theta}(T_k^{N,-})- \Lambda^-_{\theta}(T).$$ 
The terms $\ell_T^+(\theta)$ and $\ell_T^-(\theta)$ represent the log-likelihood contributions associated with the excitatory and inhibitory jumps, respectively, with $\lambda^+$ and $\lambda^-$ denoting their corresponding intensity functions, and $\Lambda_\theta^+(T)$ and $\Lambda_\theta^-(T)$ their respective compensators. The estimator is defined as the minimiser of the negative log-likelihood given in~\eqref{MLE}. Since the log-likelihood function may exhibit multiple local maxima, gradient-based algorithms are not well suited for maximum likelihood estimation. We therefore employ the derivative-free Nelder–Mead simplex algorithm (NM).
\begin{table}[htbp]
\centering
\caption{MSEs of the proposed and MLE estimators for different values of $T$, with $N=2000$ and $M=500$ Monte Carlo replications.}
\label{tab:MSE_comparison}
\renewcommand{\arraystretch}{1.5}
\setlength{\tabcolsep}{10pt}
\resizebox{\textwidth}{!}
{\begin{tabular}{c cc cc cc}
\toprule
$T$
& $\widehat{v}_0$
& $\widehat{v}_{0,MLE}$
& $\widehat{\alpha}$
& $\widehat{\alpha}_{MLE}$
& $\widehat{a}$
& $\widehat{a}_{MLE}$ \\
\midrule
50  &$\texttt{5.8198e-2}$ &$\texttt{1.1006e-4}$ &$\texttt{1.0460e-1}$ & $ \texttt{7.1778e-5}$&$ \texttt{1.8365e-6}$ &$ \texttt{1.5971e-5}$ \\
100 &$\texttt{3.0804e-2}$ & $\texttt{5.9176e-5}$& $\texttt{4.8460e-2}$&$\texttt{3.4217e-5}$ &$\texttt{1.0113e-6}$  &$\texttt{9.1004e-6}$ \\
200 & $\texttt{1.5552e-2}$&$\texttt{3.1901 e-5}$ & $\texttt{2.2880 e-2}$&$\texttt{1.7493e-5}$&$\texttt{4.7306e-7}$ &$\texttt{4.8867e-6}$  \\
300 & $\texttt{1.0085e-2}$&$\texttt{2.1116e-5}$& $\texttt{1.4783e-2} $&$\texttt{1.1877e-5}$& $\texttt{3.2371e-7}$& $\texttt{3.2442e-6}$\\
\bottomrule
\end{tabular}}
\end{table}

\begin{table}[htbp]
\centering
\caption{MSEs of the proposed and MLE estimators for different values of $N$, with $T=300$ and $M=500$ Monte Carlo replications.}
\label{tab:MSE_N}
\renewcommand{\arraystretch}{1.5}
\setlength{\tabcolsep}{10pt}
\resizebox{\textwidth}{!}
{\begin{tabular}{c cc cc cc}
\toprule
$N$
& $\widehat{v}_0$
& $\widehat{v}_{0,MLE}$
& $\widehat{\alpha}$
& $\widehat{\alpha}_{MLE}$
& $\widehat{a}$
& $\widehat{a}_{MLE}$ \\
\midrule
50   &$\texttt{1.1129e-2}$  &  $\texttt{7.0112e-4}$ &$\texttt{1.4297e-2}$ & $\texttt{4.4625e-4}$ &$\texttt{1.2440e-5}$ &$\texttt{1.0904e-4}$  \\
100  &$\texttt{1.0376e-2}$  &$\texttt{4.2307e-4}$  &$\texttt{1.4034e-2}$ &  $\texttt{2.4086e-4}$& $\texttt{6.4735e-6}$ & $\texttt{6.2493e-5}$ \\
500  & $\texttt{9.4311e-3}$ & $\texttt{8.3054e-5}$& $\texttt{1.3683e-2}$& $\texttt{4.4431e-5}$&$\texttt{1.1347e-6}$  &$\texttt{1.2983e-5}$ \\
1000 &$\texttt{1.0119e-2}$ & $\texttt{3.8011e-5}$&$\texttt{1.4010e-2}$ &$\texttt{1.9507e-5}$ &$\texttt{6.3507e-7}$ &$\texttt{5.8874e-6}$ \\
2000 &$\texttt{1.0085e-2}$&$\texttt{2.1116e-5}$& $\texttt{1.4783e-2} $&$\texttt{1.1877e-5}$& $\texttt{3.2371e-7}$& $\texttt{3.2442e-6}$ \\
\bottomrule
\end{tabular}}
\end{table}
Table \ref{tab:MSE_comparison} presents the mean squared errors (MSEs) of the proposed estimator and the MLE across different observation horizons $T,$ with $N=2000$ and $M=500$ Monte Carlo replications. Both estimation procedures exhibit a decreasing MSE as $T$ increases. The MLE has smaller MSEs for the estimators of  $v_0$ and $\alpha$ while the proposed estimator has a smaller MSE for $a.$
Table \ref{tab:MSE_N} shows the MSEs of the two estimators for different values of $N$, with $T=300$ and $M=500$ Monte Carlo replications. The MSE of the MLE decreases as $N$ increases for all parameters. For the proposed estimator, the MSE decreases for~$a$ as~$N$ increases  but remains unchanged for~$v_0$ and~$\alpha$ (see Figures~\ref{MSE4:graph} and~\ref{MSE5:graph}).

\subsection{Discussion}
Although the MLE may provide more accurate estimates for $\alpha,v_0$, its implementation becomes challenging in high-dimensional settings. In such settings, the process typically exhibits a large number of jumps, and its trajectory becomes chaotic and diffusive, with behavior close to that of a CIR process. Consequently, identifying the jumps and extracting information about their occurrence times and types from the observed data become difficult. The proposed estimators do not require the jump times, unlike the MLE which relies on them. In practice, the jump times can be obtained by observing all event trains in the network, but this is costly in terms of both data collection and computational resources, making the MLE less practical in high-dimensional settings. Moreover, in the existing literature, MLE methods generally rely on numerical optimization procedures, whereas our estimators are theoretically justified and admit explicit expressions that are easy to implement. More generally, in high-dimensional multivariate Hawkes models, the large number of interaction parameters and the poor convexity properties of the likelihood function can make numerical optimization challenging. 

\section{Proof of Theorem~\ref{T1}}\label{sec:31}

One of the major difficulty in our approach is that the square-root function appearing in the SDE of the CIR~\eqref{eq:Xbar} is not Lipschitz continuous at zero. To handle this, we replace the coefficients of the corresponding SDEs by suitably truncated versions so that both the derivatives of the coefficients are bounded. For this purpose, $\text{let}\, \delta_{\varepsilon}:\mathbb{R}\to\mathbb{R}
\text{ be a } C^\infty \text{, non-decreasing function such that }$
$$\delta_{\varepsilon}(x)=
\begin{cases}
x, & x\geq \varepsilon,\\
\dfrac{3\varepsilon}{4}, & x\leq \dfrac{\varepsilon}{2}.
\end{cases}$$
and such that, there exists a constant \(M>0\) satisfying
$$
\sum_{k=1}^{3}\left\|\delta_{\varepsilon}^{(k)}\right\|_{\infty}\leq M.$$
For each $\varepsilon>0$, we define the processes $X^{N,\varepsilon}$ and $\bar{X}^{\varepsilon}$ in the same way as $X^N$ and $\bar{X}$, respectively, by replacing the jump rate functions by $f_\eps = a\,\delta_\eps$ and the diffusive coefficient by $\sqrt{2a \delta_\eps}$. 
 More explicitly, $X^{N,\varepsilon}$ is defined as the solution of
$$dX^{N,\varepsilon}_t=-\alpha(X^{N,\varepsilon}_t-v_0)dt+\frac{1}{\sqrt{N}}\sum_{j=1}^N\int_{\mathbb{R}_+}\mathbf{1}_{\{z\leq f_\varepsilon(X^{N,\varepsilon}_{t-})\}}d\pi_+^j(dt,dz)-\frac{1}{\sqrt{N}}\sum_{j=1}^N\int_{\mathbb{R}_+}\mathbf{1}_{\{z\leq f_\varepsilon(X^{N,\varepsilon}_{t-})\}}d\pi_-^j(dt,dz)$$
and 
$$d\bar{X}^{\varepsilon}_t=-\alpha (\bar{X}^{\varepsilon}_t-v_0)\,dt
+ \sqrt{2f_{\varepsilon}(\bar{X}^{\varepsilon}_t)}\,dB_t
= -\alpha (\bar{X}^{\varepsilon}_t-v_0)\,dt
+ \sqrt{2a\delta_{\varepsilon}(\bar{X}^{\varepsilon}_t)}\,dB_t.$$

Another regularization that we use concern the processes $I^N_t,J^N_t,\bar I_t,\bar J_t.$ The factor $1/t$ in the expressions of these processes creates an irregularity in the non-homogeneous dynamics at $t=0$. To handle it, let us introduce another parameter~$\gamma>0$ and define a "regular" version of these processes:
\begin{align*}
I^{N,\gamma}_t=& \frac{1}{t+\gamma} \int_{0}^{t}  X^{N}_sds,~~J^{N,\gamma}_t= \frac{1}{t+\gamma} \int_{0}^{t}  \frac{1}{X^{N}_s+\eta_{N}} \, ds,\\
\bar{I}^{\gamma}_t=& \frac{1}{t+\gamma} \int_{0}^{t} \bar X_sds,~~
\bar{J}^{\gamma}_t= \frac{1}{t+\gamma} \int_{0}^{t}\frac{1}{\bar X_s+\eta_N}ds.
\end{align*}

In addition, we need to introduce a version of these processes with both the regularization in~$\eps$ and in~$\gamma$.
\begin{align*}
I^{N,\gamma,\eps}_t=& \frac{1}{t+\gamma} \int_{0}^{t}  X^{N,\eps}_sds,~~J^{N,\gamma,\eps}_t= \frac{1}{t+\gamma} \int_{0}^{t}  \frac{1}{\delta_\eps(X^{N,\eps}_s)+\eta_{N}} \, ds,\\
\bar{I}^{\gamma,\eps}_t=& \frac{1}{t+\gamma} \int_{0}^{t} \bar X^{\eps}_sds,~~
\bar{J}^{\gamma,\eps}_t= \frac{1}{t+\gamma} \int_{0}^{t}\frac{1}{\delta_\eps(\bar X_s)+\eta_N}ds.
\end{align*}

Roughly speaking, the processes above are defined replacing the processes $X^N,\bar X$ by $X^{N,\eps},\bar X^\eps$ with two particularities:
\begin{itemize}
    \item in the expressions of $J^{N,\gamma,\eps},\bar J^{\gamma,\eps}$, we re-use the function $\delta_\eps$ to guarantee the process not to explode (for $\bar J,J^N$ Lemma~\ref{L1}.$(iii)-(iv)$ guarantee that the finiteness of the negative moments of $\bar X$ and $X^N+\eta_N$, but the regularized versions~$\bar X^\eps,X^{N,\eps}$ can reach zero),
    \item in the expression of $\bar J^{\gamma,\eps},\bar J^\gamma,$ we add a dependency w.r.t.~$\eta_N$ to simplify some technical argument in the proofs. We do not emphasize this dependency in the notation to simplify.
\end{itemize}

Recalling that Theorem~\ref{T1} states a control between the distributions of $(X^N_T,I^N_T,J^N_T)$ and $(\bar X_T,\bar I_T,\bar J_T)$, its proof is divided in five steps, consisting in controlling successively and pair-wisely:


\[
\renewcommand{\arraystretch}{1.5}
\begin{array}{c@{\qquad}c@{\qquad}c@{\qquad}c@{\qquad}c}
(X_T^N,I_T^N,J_T^N)
& \xrightarrow{\phantom{\qquad\qquad}}
& (X_T^N,I_T^{N,\gamma},J_T^{N,\gamma})
& \xrightarrow{\phantom{\qquad\qquad}}
& (X_T^{N,\varepsilon},
   I_T^{N,\gamma,\varepsilon},
   J_T^{N,\gamma,\varepsilon})
\\[4ex]
&&&& \Big\downarrow
\\[4ex]
(\bar X_T,\bar I_T,\bar J_T)
& \xleftarrow{\phantom{\qquad\qquad}}
& (\bar X_T,\bar I_T^\gamma,\bar J_T^\gamma)
& \xleftarrow{\phantom{\qquad\qquad}}
& (\bar X_T^\varepsilon,
   \bar I_T^{\gamma,\varepsilon},
   \bar J_T^{\gamma,\varepsilon})
\end{array}
\]

These controls are formally stated in the next three propositions. Notice that, once again, it is crucial to have an explicit convergence speed for each control.

\begin{proposition}\label{P1}
There exist~$C,C_2>0$ such that, for any test-function $g:\mathbb{R}^3\to\r$ of class $C^3_b$, for all~$T,\eps,\gamma,x>0,$ and $N\in\mathbb{N}^*$
$$
\left|  \mathbb{E}_x\left[
g\left(X_T^{N,\varepsilon}, I_T^{N,\gamma,\varepsilon}, J_T^{N,\gamma,\varepsilon}\right)
\right]
-
\mathbb{E}_x\left[
g\left(\overline{X}_T^{\varepsilon}, \overline{I}_T^{\gamma,\varepsilon}, \overline{J}_T^{\gamma,\varepsilon}\right)
\right]\right|\leq \frac{C}{\sqrt{N}}||g||_{3,\infty} e^{C_2\,T/\eps}(1+x). 
$$
\end{proposition}
\begin{proposition}\label{P2}
There exists~$C>0$ such that, for any continuous and bounded test-function $g:\mathbb{R}^3\to \mathbb{R}$, for all~$T,\eps,\gamma,x>0,$ and $N\in\mathbb{N}^*$
\begin{align*}
&\left| 
\mathbb{E}_x\left[g\left(\bar X_T^{\varepsilon},\bar I_T^{\gamma,\varepsilon}, \bar J_T^{\gamma,\varepsilon}\right)\right]-\mathbb{E}_x\left[g\left(\bar X_T, \bar I_T^{\gamma}, \bar J_T^{\gamma}\right)\right]\right|\leq C \varepsilon\|g\|_{\infty} \left[\frac{1}{x}+ T\left(1+\frac{1}{x}\right) +\sqrt{T}\left(1+\frac{1}{x^{3/2}}\right)\right],\\
&\left| 
\mathbb{E}_x\left[g\left(X_T^{N,\varepsilon},I_T^{N,\gamma,\varepsilon},J_T^{N,\gamma,\varepsilon}\right)\right]-\mathbb{E}_x\left[g\left(X_T^{N},I_T^{N,\gamma},J_T^{N,\gamma}\right)\right]\right|\\
&\hspace*{3.9cm} \leq
C \varepsilon\|g\|_{\infty} \left[\frac{1}{x+\eta_N}+ T\left(1+\frac{1}{x+\eta_N}\right) +\sqrt{T}\left(1+\frac{1}{(x+\eta_N)^{3/2}}\right)\right].
\end{align*}
\end{proposition}
\begin{proposition}\label{P3}
There exists~$C>0$ such that for all~$T,\gamma,x>0,$ and $N\in\mathbb{N}^*$,
the following $L^1$ controls hold
\begin{align*}
    &\mathbb{E}_x\left[\left|\left|\left(\bar X_T, \bar I^{\gamma}_T, \bar J^{\gamma}_T\right) - \left(\bar X_T, \bar I_T, \bar J_T\right)\right|\right|_1\right]\leq C\frac{\gamma}{T}\left(1 + x + \frac1x\right) + C\eta_N \left(1+\frac{1}{x^2}\right),\\
    &\mathbb{E}_x\left[\left|\left|\left(X^{N}_T,I^{N,\gamma}_T,J^{N,\gamma}_T\right) - \left(X^{N}_T,I^{N}_T,J^{N}_T\right)\right|\right|_1\right]\leq C\frac{\gamma}{T}\left(1 + x + \frac{1}{x+\eta_N}\right).
\end{align*}
\end{proposition}

Throughout the proofs of these propositions, it is required to have suitable controls on our processes~$X^N,\bar X$. Formally, some a priori estimates are given in the lemma below. The proof is postponed to the Appendix.

\begin{lemma}\label{L1}$ $
\begin{itemize}
\item[$(i)$] For any~$p,t>0$, there exists~$C_{p,t}>0,$ such that for all~$x>0,$
$$\mathbb{E}_x\left[|\bar X_t|^{p}\right]\leq C_{p,t}\left(1+x^p\right)\textrm{ and }\underset{N\in\mathbb{N}^*}{\sup}\mathbb{E}_x\left[\left|X^N_t\right|^{p}\right]\leq C_{p,t}\left(1+x^p\right).$$
\item[$(ii)$] For any~$x,t,N,$
$$\mathbb{E}_x\left[X^N_t\right] = \mathbb{E}_x\left[\bar X_t\right] = v_0 + (x-v_0)e^{-\alpha t}.$$
\item[$(iii)$] Let $p>0$, if $a\,p < \alpha\, v_0,$ then there exists $C_p>0$, such that, for all~$x>0,$
$$\underset{t\geq 0}{\sup}~\mathbb{E}_x\left[\left(\bar X_t\right)^{-p}\right] \leq C_p\left(1+x^{-p}\right).$$
\item[$(iv)$] Let $p>0$, if $a(p+1) < \alpha\, v_0,$ then there exists $C_p>0$, such that, for all~$x,t>0,$
$$\underset{t\geq 0}{\sup}~\mathbb{E}_x\left[\left(X^N_t + \eta_N\right)^{-p}\right] \leq \left(x + \eta_N\right)^{-p} +C_p\, t,$$
where we recall that $\eta_N$ can be any sequence satisfying $N^{-1/2}<<\eta_N<<1$.
\end{itemize}
\end{lemma}

Item~$(i)$ is used implicitly throughout the paper to guarantee that all the local martingales (i.e. integrals w.r.t. Brownian motions or compensated Poisson measures) are martingales. Item~$(ii)$ is useful to have finer control on the first moments of the processes. And the negative moments given in~$(iii)$ and~$(iv)$ are mainly used in the proof of Proposition~\ref{P2} to control the probability that the processes come close to zero. They are also used in Proposition~\ref{P3} to control the moments of~$J^N_T$ and~$\bar J_T$.

\subsection{Proof of Proposition~\ref{P1}} 
In this part, we will work with the following processes 
$\Bar{Y}^{\varepsilon}=\bar{X}^{\varepsilon}-v_0,$ and $Y^{N,{\varepsilon}}=X^{N,{\varepsilon}}-v_0,$ so that 
$$d\bar{Y}^{\varepsilon}_t
= -\alpha \bar{Y}_t\,dt
+ \sqrt{2a\,\delta_{\varepsilon}(\bar{Y}_t+v_0)}\,dB_t,$$
and
$$dY^{N,\varepsilon}_t = -\alpha Y^N_t dt + \frac1{\sqrt{N}}\sum_{j=1}^N \int_{\r_+}\uno{z\leq f_{\varepsilon}(Y^N_{t-}+v_0)}d\pi_+^j(dt,dz) - \frac1{\sqrt{N}}\sum_{j=1}^N \int_{\r_+}\uno{z\leq f_{\varepsilon}(Y^N_{t-}+v_0)}d\pi_-^j(dt,dz).$$
where $f_{\varepsilon}(x) = a\delta_\eps(x).$ 
 For the sake of notation, let us introduce $\xi^N := (Y^{N,\eps},I^{N,\gamma,\eps},J^{N,\gamma,\eps})$ and $\bar \xi := (\bar Y^{\eps},\bar I^{\gamma,\eps},\bar J^{\gamma,\eps}),$ which will be used throughout the proof.
\begin{lemma}
Let $A^{N}$ be the generator of $\xi^N$ and $\Bar{A}$ the generator of $\bar \xi.$ For all \(g \in C_b^3(\mathbb{R})\), and for any~$t>0$ and $x,y,z>0$, 
$$\left|A_t^N g(x,y,z)-\bar{A}_tg(x,y,z)\right|\leq \frac{f_{\varepsilon}(x+v_0)}{3\sqrt{N}}\left\|\partial_1^3 g\right\|_\infty.$$   
\end{lemma}
\begin{proof}
To start, we have $$dI_t^{N,\gamma,\eps}=\left[\frac{Y^{N,\varepsilon}_t+v_0}{t+\gamma}-\frac{I^{N,\gamma,\eps}_t}{t+\gamma}\right]dt~~\textrm{ and }~~dJ^{N,\gamma,\eps}_t=\left(\frac{1}{(t+\gamma)(\delta_\eps(Y^{N,\varepsilon}_t+v_0)+\eta_{N})}
-\frac{J^{N,\gamma,\eps}_t}{t+\gamma}\right)dt.$$
Let $g \in C_b^{3}(\mathbb{R})$ be a test function, by Ito's formula, we have 
\begin{align*}
g(\xi^N_t)
=&a{}g(\xi^N_0)
-\alpha\int_0^t Y_s^{N,\varepsilon}\partial_1g(\xi^N_s)ds +
\int_0^t\frac{Y_s^{N,\varepsilon}+v_0-I_s^{N,\gamma,\eps}}{s+\gamma}\partial_2g(\xi^N_s)ds\\&+\int_0^t\left(\frac{1}{(t+\gamma)(\delta_\eps(Y^{N,\varepsilon}_t+v_0)+\eta_{N})}
-\frac{J^{N,\gamma,\eps}_t}{t+\gamma}\right)\partial_3g(\xi^N_s)ds\\&+ \sum_{j=1}^N \int_{0}^{t} f_{\varepsilon}\left(Y^{N,\varepsilon}_{s^-}+v_0\right) \left[g\left(Y^{N,\varepsilon}_{s^-}+\frac{1}{\sqrt{N}},I^{N,\gamma,\eps}_{s^-},J^{N,\gamma,\eps}_{s^-}\right) -g\left(\xi^N_{s^-}\right)\right]ds \\ &+\sum_{j=1}^N \int_{0}^{t} f_{\varepsilon}\left(Y^{N,\varepsilon}_{s^-}+v_0\right) \left[g\left(Y^{N,\varepsilon}_{s^-}-\frac{1}{\sqrt{N}},I^{N,\gamma,\eps}_{s^-},J^{N,\gamma,\eps}_{s^-}\right) -g\left(\xi^N_{s^-}\right)\right]ds +M^N_t,
\end{align*}
where $M^N$ is a local martingale. Then, by the martingale characterization of the generator,
\begin{align*}
    A^N_t g(x,y,z)=&  -\alpha x \partial_1 g(x,y,z) +\left[\frac{x+v_0-y}{t+\gamma}\right] \partial_2 g(x,y,z)\\& +\left[\frac{1}{(t+\gamma)(\delta_\eps(x+v_0)+\eta_N)}-\frac{z}{t+\gamma}\right] \partial_3 g(x,y,z)\\&+  N
f_{\varepsilon}\left(x+v_0\right)
\left[g\left(x + \frac{1}{\sqrt{N}}, y,z\right)-g\left(x,y,z\right)
\right]\\&+N
f_{\varepsilon}\left(x+v_0\right)
\left[g\left(x - \frac{1}{\sqrt{N}}, y,z\right)-g\left(x, y,z\right)
\right]
\end{align*}

Using again Ito’s formula and the martingale characterization of the generator to calculate $\Bar{A}$ the generator of $\bar \xi := (\Bar{Y}^{\varepsilon},\Bar{I}^{\gamma,\eps},\Bar{J}^{\gamma,\eps}),$ we obtain
\begin{align*}
 \Bar{A}_tg(x,y,z)=&-\alpha x \partial_1 g(x,y,z)+\Big[\frac{x+v_0-y}{t+\gamma}\Big]\partial_2 g(x,y,z)\\
 &+\Big[\frac{1}{(t+\gamma)(\delta_\eps(x+v_0)+\eta_N)}-\frac{z}{t+\gamma}\Big]\partial_3 g(x,y,z)
 +\partial_1^2g(x,y,z)f_{\varepsilon}(x+v_0).
\end{align*}
Hence
\begin{align*}
\left|A_t^N g(x,y,z)-\bar{A}_tg(x,y,z)\right|
&=
f_{\varepsilon}(x+v_0)\Bigg|
N\left[
g\left(x+\frac{1}{\sqrt{N}},y,z\right)-g(x,y,z)
\right]  \\
&\qquad\qquad
+N\left[
g\left(x-\frac{1}{\sqrt{N}},y,z\right)-g(x,y,z)
\right]
-\partial_1^2 g(x,y,z)
\Bigg|.
\end{align*}
From Taylor-Lagrange’s inequality, we get
\begin{align*}
&\left|A_t^N g(x,y,z)-\bar{A}_t g(x,y,z)\right|
\\
&\qquad\leq f_{\varepsilon}(x+v_0)\Bigg|
N\Bigg[
\frac{1}{\sqrt{N}}\partial_1 g(x,y,z)
+\frac{1}{2N}\partial_1^2g(x,y,z)
+\frac{1}{6N^{3/2}}\left\|\partial_1^3 g\right\|_\infty
\Bigg] \\
&\qquad\quad
+N\Bigg[
-\frac{1}{\sqrt{N}}\partial_1 g(x,y,z)
+\frac{1}{2N}\partial_1^2g(x,y,z)
+\frac{1}{6N^{3/2}}\left\|\partial_1^3 g\right\|_\infty
\Bigg] 
-\partial_1^2 g(x,y,z)
\Bigg| \\
&\qquad\leq \frac{f_{\varepsilon}(x+v_0)}{3\sqrt{N}}\left\|\partial_1^3 g\right\|_\infty,
\end{align*}
which is the desired control.
\end{proof}

\begin{lemma}\label{lemma:semigr}
Let $(P^N_{r,t})_{r\leq t}$ and $(\bar P_{r,t})_{r\leq t}$
denote the (non-homogeneous) Markov semigroups associated respectively to $\xi^N$ and $\bar \xi.$ Then, there exists some~$C,C_2>0$, such that for all~$g\in C^3_b(\r),N\in\mathbb{N}^*,0\leq r\leq t,\eps>0,\gamma>0,(x,y,z)\in(0,+\infty),$ 
\begin{equation*}
\left|\bar{P}_{r,t}g(x,y,z)-P_{r,t}^N g(x,y,z)\right|\leq \frac{3}{\sqrt{N}} C ||g||_{3,\infty} e^{C_2 t/\eps}.
\end{equation*}
\end{lemma}

\begin{proof}
For the convergence of the semigroups, we use Lemma~\ref{lem:trotter}:
\begin{equation*}
    \left(\bar{P}_{(r,t)} - P_{(r,t)}^{N,}\right)g(x,y,z)=\int_0^tP_{r,{t-s}}^{N}(\bar{A}_r-A^{N}_r)\Bar{P}_{(r,s)}g(x,y,z)ds,
\end{equation*}
Using lemma \ref{lem:diffsemi},
\begin{align*}
\left|\bar{P}_{r,t}g(x,y,z)-P_{r,t}^{N}g(x,y,z)\right|
&=\Bigg|\int_0^t P_{r,t-s}^{N}\left(\bar{A}_r-A_r^{N}\right)\bar{P}_{r,s}g(x,y,z)\,ds\Bigg|\\
&\leq\int_0^t\mathbb{E}_{(x,y,z)}\Bigg[\Bigg|\bar{A}_r\left(\bar{P}_{r,s}g\right)\left(\xi^N_{t-s}\right)-A_r^{N}\left(\bar{P}_{r,s}g\right)\left(\xi^{N}_{t-s}\right)\Bigg|\Bigg]\,ds\\
&\leq \frac{3}{\sqrt{N}} \int_0^t\mathbb{E}_{(x,y,z)} [f_{\varepsilon}(Y_{t-s}^{N,\varepsilon}+v_0)] \left\|\partial_x^3 \bar{P}_{r,s}^{\varepsilon}\right\|_\infty ds\\
&\leq\frac{3}{\sqrt{N}} C_1 (1+x) ||g||_{3,\infty}\int_0^t e^{C_2 s/\eps}ds\\
& \leq\frac{3}{\sqrt{N}} C (1+x) ||g||_{3,\infty} e^{C_2 t/\eps},
\end{align*}
where, to obtain the second last inequality above, we have used that
$$\esp{f_\eps (Y^{N,\eps}_t)} \leq \esp{Y^{N,\eps}_t} = v_0 + (\esp{Y_0}-v_0)e^{-\alpha t} \leq C(1+\esp{Y_0}).$$
This proves the lemma.
\end{proof}

\subsection{Proof of Proposition~\ref{P2}}

{\it Step~1.} Let us prove that, there exists~$C>0$ such that for all~$t,x>0,$
\begin{equation}\label{eq:step1}
    \mathbb{E}_x\left[\underset{s\leq t}{\sup}\,\,\frac1{\bar X_s}\right] \leq C\left[\frac1x + t(1+1/x) + \sqrt{t(1+1/x^3)}\right].
\end{equation}

Applying Itô's formula, we obtain
\begin{align*}
\frac{1}{\bar X_t}&= \frac{1}{x}-\int_0^t \frac{1}{{\bar X_s}^2}d\bar X_s+\int_0^t \frac{2a}{{\bar X_s}^2} ds\\&= \frac{1}{x}+\int_0^t\left[ \frac{\alpha}{\bar X_s} +\frac{2a-\alpha v_0}{{\bar X_s}^2}\right]ds -\int_0^t \frac{\sqrt{2a}}{{\bar X}^{3/2}} dB_s\\&= \frac{1}{x}+\int_0^t\frac{\alpha}{\bar X_s}ds -\int_0^t \frac{\sqrt{2a}}{{\bar X}^{3/2}} dB_s
\end{align*}
Therefore, applying Fubini's theorem followed by the B\"urkholder-Davis-Gundy's inequality, and then Jensen's inequality, we obtain
\begin{align*}
   \mathbb{E}\left[\sup_{s\leq t}\frac{1}{\bar X_s}\right]&\leq\mathbb{E}\left[\frac{1}{x}+\alpha\int_0^t\frac{1}{\bar X_s} ds+\sup_{s\leq t}|M_s|\right]\\&\leq\frac{1}{x}+\alpha\int_0^t
\mathbb{E}\left[\frac{1}{\bar X_s}\right] ds +\mathbb{E}\left[\sup_{s\leq t}|M_s|\right]\\&\leq\frac{1}{x}+\alpha\int_0^t
\mathbb{E}\left[\frac{1}{\bar X_s}\right] ds +C'\mathbb{E}\left[\left(\int_0^t \frac{2a}{\bar X_s^3} ds
\right)^{1/2}\right]\\&\leq\frac{1}{x}+\alpha\int_0^t
\mathbb{E}\left[\frac{1}{\bar X_s}\right] ds+ C'\sqrt{2a}\left(\int_0^t\mathbb{E}\left[\frac{1}{\bar X_s^3}\right]ds\right)^{1/2}.
\end{align*}
Then, Lemma~\ref{L1}.$(iii)$ entails~\eqref{eq:step1}.

{\it Step~2.} Now, let us prove the first inequality of Proposition~\ref{P2}. Let $$A_{\varepsilon}=\left\{\inf_{s\leq T}\bar X_s>\varepsilon\right\},$$
and notice that, on the event~$A_\eps$, the processes $\left(\bar X_s\right)_{s\leq T}$ and $\left(\bar X^\eps_s\right)_{s\leq T}$
are equal. Consequently
\begin{align*}
&\left| 
\mathbb{E}\left[g\left(\bar X_T^{\varepsilon},\bar I_T^{\gamma,\varepsilon}, \bar J_T^{\gamma,\varepsilon}\right)\right]-\mathbb{E}\left[g\left(\bar X_T, \bar I_T^{\gamma}, \bar J_T^{\gamma}\right)\right]\right|\\
&=
\left|\mathbb{E}\left[g\left(\bar X_T^{\varepsilon},\bar I_T^{\gamma,\varepsilon}, \bar J_T^{\gamma,\varepsilon}\right)\mathbf{1}_{A_{\varepsilon}}\right]
+\mathbb{E}\left[g\left(\bar X_T^{\varepsilon},\bar I_T^{\gamma,\varepsilon}, \bar J_T^{\gamma,\varepsilon}\right)\mathbf{1}_{A_{\varepsilon}^c}\right]\right.\\
&~~~~~~\left.- \mathbb{E}\left[g\left(\bar X_T, \bar I_T^{\gamma}, \bar J_T^{\gamma}\right)\mathbf{1}_{A^{\varepsilon}}\right]-\mathbb{E}\left[g\left(\bar X_T, \bar I_T^{\gamma}, \bar J_T^{\gamma}\right)\mathbf{1}_{A_{\varepsilon}^c}\right]\right|\\
&= \left|\mathbb{E}\left[g\left(\bar X_T^{\varepsilon},\bar I_T^{\gamma,\varepsilon}, \bar J_T^{\gamma,\varepsilon}\right)\mathbf{1}_{A_{\varepsilon}^c}\right] - \mathbb{E}\left[g\left(\bar X_T, \bar I_T^{\gamma}, \bar J_T^{\gamma}\right)\mathbf{1}_{A_{\varepsilon}^c}\right]\right|\\
&\leq \left|\mathbb{E}\left[g\left(\bar X_T^{\varepsilon},\bar I_T^{\gamma,\varepsilon}, \bar J_T^{\gamma,\varepsilon}\right)\mathbf{1}_{A_{\varepsilon}^c}\right]\right|+\left|\mathbb{E}\left[g\left(\bar X_T, \bar I_T^{\gamma}, \bar J_T^{\gamma}\right)\mathbf{1}_{A_{\varepsilon}^c}\right]\right|\\
&\leq 2\|g\|_{\infty}\mathbb{P}\left(\inf_{s\leq T}\bar X_s\leq\varepsilon\right)
\end{align*}
Using Markov's inequality, we obtain
\begin{align*}
\mathbb{P}\left(\inf_{s\leq T}\bar X_s\leq\varepsilon\right)&\leq \mathbb{P}\left(\sup_{s\leq T}\,\,\frac{1}{\bar X_s}\geq \frac{1}{\varepsilon}\right)\leq \varepsilon\,\mathbb{E}\left[\sup_{s\leq T}\,\,\frac{1}{\bar X_s}\right]
\end{align*}
and~\eqref{eq:step1} proved at {\it Step~1} allows to conclude.

{\it Step~3.} In this last step, we prove the second inequality of the statement of Proposition~\ref{P2}. As in {\it Step~1}, using Lemma~\ref{L1}.$(iv)$ instead of Lemma~\ref{L1}.$(iii)$, we can prove that
$$\mathbb{E}_x\left[\underset{s\leq T}{\sup}\,\,\frac1{X^N_s + \eta_N}\right] \leq C\left[\frac1{x+\eta_N} + T\left(1+1/(x+\eta_N)\right) + \sqrt{T\left(1+1/(x+\eta_N)^3\right))}\right].$$

Then, the same calculation as the ones at {\it Step~2} allow to conclude the proof using the inequality above.

\subsection{Proof of Proposition~\ref{P3}}


Relying on Lemma~\ref{L1} to control the moments of $\bar X_T$ and $1/\bar X_T$, we have
  \begin{align*}
      &\mathbb{E}\left[\left|\left|\left(\bar X_T, \bar I^{\gamma}_T, \bar J^{\gamma}_T\right) - \left(\bar X_T, \bar I_T, \bar J_T\right)\right|\right|_1\right]\\
      &\leq \esp{\left|\frac1T\int_0^T \bar X_sds - \frac1{T+\gamma}\int_0^T \bar X_sds\right|} + \esp{\left|\frac1T\int_0^T \frac1{\bar X_s}ds - \frac1{T+\gamma}\int_0^T \frac1{\bar X_s+\eta_N}ds\right|}\\
      &\leq \left|\frac1T - \frac{1}{T+\gamma}\right|\esp{\int_0^T \bar X_sds} + \left|\frac1T - \frac{1}{T+\gamma}\right|\esp{\int_0^T \frac{1}{\bar X_s}ds} + \frac{1}{T+\gamma}\esp{\int_0^T\left|\frac{1}{\bar X_s}-\frac{1}{\bar X_s + \eta_N}\right|ds}\\
      &\leq C\frac{\gamma}{T}\left(1+x+\frac1x\right) + \frac1{T+\gamma} \int_0^T \esp{\frac{\eta_N}{(\bar X_s)^2}}ds\\
      &\leq C\frac{\gamma}{T}\left(1+x+\frac1x\right) + C\eta_N\left(1+1/x^2\right).
      \end{align*}
      
And, similarly,
\begin{align*}
&\mathbb{E}\left[\left|\left|\left(X^N_T, I^{N,\gamma}_T, J^{N,\gamma}_T\right) - \left(X^N_T, I^N_T, J^N_T\right)\right|\right|_1\right]\\
&\leq \esp{\left|\frac1T\int_0^T X^N_sds - \frac1{T+\gamma}\int_0^T X^N_sds\right|} + \esp{\left|\frac1T\int_0^T \frac1{X^N_s+\eta_N}ds - \frac1{T+\gamma}\int_0^T \frac1{\bar X^N_s+\eta_N}ds\right|}\\
&\leq \left|\frac1T - \frac{1}{T+\gamma}\right|\esp{\int_0^T X^N_sds} + \left|\frac1T - \frac{1}{T+\gamma}\right|\esp{\int_0^T \frac{1}{X^N_s+\eta_N}ds}\\
&\leq C \frac{\gamma}{T} \left(1 + x + 1/(x+\eta_N)\right).
\end{align*}

\subsection{Proof of Theorem~\ref{T1}}
\label{sec:thmt1}

Combining the results of the Propositions~\ref{P1},~\ref{P2} and~\ref{P3}, we have that, for any test-function $g\in C^3_b(\r^3)$ and any $x>0,$ there exists $C_{x,g}>0$ such that, for all~$T,\gamma,\eps>0,$
\begin{equation}\label{eq:thm21}
\left|\mathbb{E}_x\left[g\left(X_T^N,I_T^N,J_T^N\right)\right] - \mathbb{E}_x\left[g\left(\bar X_T,\bar I_T,\bar J_T\right)\right]\right| \leq C_{x}||g||_{3,\infty}\left(N^{-1/2} e^{C_2\,T/\eps} + T\,\eps + \frac{\gamma}{T} + \eta_N\right).
\end{equation}

In particular, no constraints come from~$\gamma,\eta_N.$ For the expression above to vanish, one needs to chose $\eps_T>0$ such that
$$\eps_T = o(1/T)$$
as $T$ goes to infinity, and then, to chose $N_T$ such that
$$N_T^{-1/2} e^{C_2\,T/\eps_T} \underset{T\to\infty}{\longrightarrow}0.$$

In particular, writing $N_T$ such that
$$N_T\geq e^{T^2 M_T},$$
with $M_T$ going to infinity as $T$ goes to infinity, and choosing $$\eps_T = \frac{1}{T\sqrt{M_T}}$$
make the RHS of~\eqref{eq:thm21} vanish.

\section{Proof of Theorem~\ref{T2}: consistencies}\label{sec:33}

Let us begin by proving that the pair $(I^N_T,J^N_T)$ converges in distribution to a deterministic value as $N,T$ tend jointly to infinity. In particular, the convergence is in probability. This result relies on two properties:
\begin{itemize}
    \item Theorem~\ref{T1}, which implies that, in distribution, the pair $(I^N_T,J^N_T)$ behaves "similarly" as the limit pair $(\bar I_T,\bar J_T)$,
    \item the convergence of the limit pair as $T$ goes to infinity.
\end{itemize}

Let us remark that, technically, Theorem~\ref{T1} is not directly a control between the distributions of $(I^N_T,J^N_T)$ and $(\bar I_T,\bar J_T)$ since we work on test-functions of class~$C^3_b$, and not on all continuous and bounded test-functions. So the next classical lemma is required to use Theorem~\ref{T1}.
\begin{lemma}\label{lem:cvLoi}
Let $d\in\mathbb{N}^*$, and $\zeta,\zeta_n$ ($n\in\mathbb{N}^*$) be r.v. taking values in~$\mathbb{R}^d$. Assume that for any $C^\infty$ compactly supported test-function~$g:\r^d\to \r,$
$$\esp{g(\zeta_n)}\underset{n\to\infty}{\longrightarrow}\esp{g(\zeta)}.$$
Then $\zeta_n$ converges to $\zeta$ in distribution as $n$ goes to infinity.
\end{lemma}

\begin{proof}[Sketch of proof]
Let $g:\mathbb{R}^d\to\mathbb{R}$ be continuous and compactly supported. It is classical that the function $g$ can be approximated uniformly by a sequence of $C^\infty_c$ functions~$(g_k)_k$, by defining $g_k$ as the convolution product between $g$ and a mollifier. Then, writing
$$\left|\esp{g(\zeta_n)}-\esp{g(\zeta)}\right| \leq \left|\esp{g(\zeta_n)}-\esp{g_k(\zeta_n)}\right|+\left|\esp{g_k(\zeta_n)}-\esp{g_k(\zeta)}\right| + \left|\esp{g_k(\zeta)}-\esp{g(\zeta)}\right|$$
implies the convergence of $\esp{g(\zeta_n)}$ to $\esp{g(\zeta)}$ for any continuous function~$g$ with a compact support.

Then, using that the vague convergence of probability measures entails classically the weak convergence, the lemma is proved.
%
\end{proof}

Now the next result comes easily from the previous results and from \cite{alaya2012parameter}.

\begin{proposition}\label{prop:cvProba}
The following convergence in probability holds true
$$\left(I^N_T,J^N_T\right)\longrightarrow \left(v_0,\frac{1}{v_0 - a/\alpha}\right)$$
as $N,T$ tend jointly to infinity, as soon as there exists $M_T\to\infty$ such that $N_T\geq e^{T^2 M_T}.$
\end{proposition}

\begin{proof}
Let us write, for any $g\in C^3_b(\mathbb{R}^2),$
\begin{multline*}
    \left|\esp{g\left(I^N_T,J^N_T\right)} - \esp{g\left(v_0,\frac{1}{v_0 - a/\alpha}\right)}\right|\\
    \leq\left|\esp{g\left(I^N_T,J^N_T\right)} - \esp{g\left(\bar I_T,\bar J_T\right)}\right|+\left|\esp{g\left(\bar I_T,\bar J_T\right)} - \esp{g\left(v_0,\frac{1}{v_0 - a/\alpha}\right)}\right|.
\end{multline*}
    
The first term of the sum above vanishes thanks to Theorem~\ref{T1}, and the second one thanks to Theorem~4.1 of \cite{alaya2012parameter}. So, we have the convergence for any test-functions in $C^3_b(\r^2)$, and so, by Lemma~\ref{lem:cvLoi}, the pair $(I^N_T,J^N_T)$ converges in distribution. Since the limit is deterministic, the convergence in probability holds true.
\end{proof}

\subsection{Consistency of $\displaystyle\widehat{v}_0$} 
Recalling that $\widehat{v}_0^N(T)=I^N_T,$ the consistency of $\widehat{v}_0^N(T)$ w.r.t.~$v_0$ is stated in Proposition~\ref{prop:cvProba}.


\subsection{Consistency of $\displaystyle\widehat{a}$}
\label{sec:consa}
 
We consider the estimator $$\widehat{a}^N(T)=\frac{[X^N]_T}{2TI_T^N}.$$
The quadratic variation of $X^N$ is characterized as
\begin{align*}
    [X^N]_t&= \sum_{0<s\leq t}(\Delta X_s^N)^2\\&=\frac{1}{N}\sum_{j=1}^{N}\int_0^t\int_{\mathbb{R}_+}\mathbf{1}_{\{\leq f(X_{s-}^N)\}}\left(d\pi_j^+(s,z)+d\pi_j^-(s,z)\right)\\&=\frac{1}{N}\sum_{j=1}^{N}\int_0^t\int_{\mathbb{R}_+}\mathbf{1}_{\{z\leq f(X_{s-}^N)\}}
\left(\widetilde{\pi}_j^+(ds,dz)+\widetilde{\pi}_j^-(ds,dz)\right)
+\frac{2}{N}\sum_{j=1}^{N}\int_0^t\int_{\mathbb{R}_+}
\mathbf{1}_{\{z\leq f(X_{s-}^N)\}}\,ds dz\\
&=M_t^N+2a\int_0^t \left(X_{s-}^N\right)_+ ds\\
&=M_t^N+2at \tilde I_t^N,
\end{align*}
with
$$\tilde I^N_t = \frac1t \int_0^t \left(X_{s}^N\right)_+ ds.$$

So,
\begin{equation}\label{eq:a0}
\widehat{a}^N(T)=\frac{[X^N]_T}{2TI_T^N}= \frac{M^N_T}{2TI^N_T}+a \frac{\tilde I^N_T}{I^N_T}.
\end{equation}
Since, by Proposition~\ref{prop:cvProba},
\begin{equation}\label{a0}
\frac{1}{I_T^N} \xrightarrow[T,N\to\infty]{\mathbb{P}} \frac{1}{v_0},
\end{equation}
it remains to prove the convergence of  $\frac{M^N_T}{T}$ to zero in $L^2$, hence in probability, and that $\tilde I^N_T/I^N_T$ converges in probability to one.
By the martingale isometry and Fubini's theorem, we have, for all~$t>0,$
\begin{equation}\label{a1}
\mathbb{E}\left[\left|\frac{M_t^N}{t}\right|^2\right]
=\frac{2a}{Nt^2}\mathbb{E}\left[\int_0^t X_{s}^N\,ds\right]=\frac{2a}{Nt^2}\int_0^t\mathbb{E}[\left(X_{s}^N\right)_+]\,ds.
\end{equation}
And, by Lemma~\ref{L1}$.(ii)$,
\begin{equation}\label{a2}
\mathbb{E}[(X^N_t)_+] \leq \mathbb{E}[X_t^N] + 1=v_0+(x_0-v_0)e^{-at}+1\leq v_0+x_0+1,
\end{equation}
where we have used that, for all~$t\geq 0,$ $X^N_t>-1/\sqrt{N},$ implying that $(X^N_t)_+\leq X^N_t+1$.

Combining \eqref{a1} and \eqref{a2}, we obtain 
$$
\mathbb{E}\left[\left|\frac{M_T^N}{T}\right|^2\right]
\leq
\frac{2a(v_0+x_0+1)}{NT}
\xrightarrow[T,N\to \infty]{}0.$$

On the other hand, denoting by
$$A^N_{\eps,T} = \left\{\underset{s\leq T}{\inf} X^N_s > \eps\right\}$$
as in the proof of Proposition~\ref{P2}, we have
$$\pro{A^N_{\eps,T}} \leq C_x\,T\,\eps,$$
and, for any~$\eta>0$,
\begin{equation}\label{eq:erra2}
\pro{\left|\frac{\tilde I^N_T}{I^N_T} - 1\right|>\eta} \leq \pro{A^N_{\eps,T}} \leq C_x\,T\,\eps\underset{(T,\eps)\to(\infty,0)}{\longrightarrow}0,
\end{equation}
when $\eps>0$ is suitably choosen depending on~$T$ as in the proof of Theorem~\ref{T1} (recalling Section~\ref{sec:thmt1}).

Therefore, from~\eqref{a1} and~\eqref{eq:erra2}, $$\widehat{a}^N(T) \xrightarrow[T,N\to\infty]{\mathbb{P}} a,$$
proving the consistency of $\widehat{a}^N(T)$ w.r.t.~$a$.

\subsection{Consistency of $\displaystyle\widehat{\alpha}$}
 
Let us recall that
$$\widehat{\alpha}^N(T) = \frac{\widehat{a}^N(T)}
{\widehat{v}_0^N(T)-1/J_T^N}.$$



By Proposition~\ref{prop:cvProba}, and by the consistencies of $\widehat{v_0}^N(T),$ $\widehat{a}^N(T)$ w.r.t.~$v_0,a$,
$$\left(\widehat{v}_0^N(T),\widehat{a}^N(T),J_T^N\right)
\xrightarrow[T,N\to\infty]{\mathbb{P}}
\left(v_0,a,\frac{1}{v_0 - a/\alpha}\right).$$
Using the Continuous Mapping Theorem, 
we deduce that
$$\widehat{\alpha}^N(T)
\xrightarrow[T,N\to\infty]{\mathbb{P}}
\frac{a}{v_0-(v_0 - a/\alpha)}=\alpha.$$

This ends the proof of Theorem~\ref{T2}.

\section{Proof of Theorem~\ref{T3}: asymptotic normalities}\label{sec:34}

As for the proof of Theorem~\ref{T2}, we begin with a result stating the joint asymptotic normality of $I^N_T,J^N_T$.

\begin{proposition}\label{prop:IJan}
The following convergence in law holds true
$$\left(\sqrt{T}\left(I^N_T-v_0\right),\sqrt{T}\left(J^N_T - \frac{1}{v_0 - a/\alpha}\right)\right)\xrightarrow{\mathcal L}  \mathcal N_2\left(0,\begin{pmatrix}
\displaystyle\frac{2av_0}{\alpha^2}
&
\displaystyle-\frac{2a}{\alpha(\alpha v_0-a)}
\\[3mm]
\displaystyle-\frac{2a}{\alpha(\alpha v_0-a)}
&
\displaystyle\frac{2a\alpha}{(\alpha v_0-a)^3}
\end{pmatrix}\right).
$$
in law as $N,T$ tend jointly to infinity, as soon as there exists $M_T\to\infty$ such that $N_T\geq e^{T^{7/2} M_T}.$
\end{proposition}

\begin{proof}
Let us first prove that 
$
\sqrt{t}\left(\bar I_t-v_0,\bar J_t-\frac{1}{v_0 - a/\alpha}\right)
$
converges in law to a normal distribution as $t$ tends to infinity with $\bar{I}_t= \frac{1}{t} \int_{0}^{t}  \bar{X}_s \, ds$ and 
$\bar{J}_t= \frac{1}{t} \int_{0}^{t} \frac{1}{ \bar{X}_s}\, ds.$
To do this, we first consider the integral form of the stochastic differential equation given by relation (\ref{eq:Xbar}) 
\[
\bar X_t
= \bar X_0
+ \alpha \int_0^t (v_0-\bar X_s)\,ds
+ \sqrt{2a}\int_0^t \sqrt{\bar X_s}\,dB_s.
\]
to deduce 
\[
\bar I_t
= v_0+\frac{\bar X_0-\bar X_t}{\alpha t}
+\frac{\sqrt{2a}}{\alpha t}
\int_0^t\sqrt{\bar X_s}\,dB_s.
\]
Second, the Itô formula with function $\ln(x)$ yields
\[
\ln(\bar X_t)-\ln(\bar X_0)
=
(\alpha v_0-a)\int_0^t\frac{1}{\bar X_s}\,ds
-\alpha t
+\sqrt{2a}\int_0^t\frac{1}{\sqrt{\bar X_s}}\,dB_s
\]
which gives us
\[
\bar J_t
=
\frac{\alpha}{\alpha v_0-a}
+
\frac{1}{\alpha v_0-a}
\frac{\ln(\bar X_t)-\ln(\bar X_0)}{t}
-
\frac{\sqrt{2a}}{\alpha v_0-a}
\frac{1}{t}
\int_0^t\frac{1}{\sqrt{\bar X_s}}\,dB_s.
\]
Hence we can write 
\[
\sqrt{t}
\begin{pmatrix}
\bar I_t-v_0\\[2mm]
\displaystyle\bar J_t-\frac{\alpha}{\alpha v_0-a}
\end{pmatrix}
=h_t+\frac{1}{\sqrt{t}}QM_t.
\]
with 
\[
M_t=
\begin{pmatrix}
\displaystyle\int_0^t\sqrt{\bar X_s}\,dB_s\\[3mm]
\displaystyle\int_0^t\frac{1}{\sqrt{\bar X_s}}\,dB_s
\end{pmatrix},
\qquad
h_t=
\begin{pmatrix}
\displaystyle\frac{\bar X_0-\bar X_t}{\alpha\sqrt{t}}\\[3mm]
\displaystyle\frac{\ln(\bar X_t)-\ln(\bar X_0)}
{(\alpha v_0-a)\sqrt{t}}
\end{pmatrix},
\qquad
Q=
\begin{pmatrix}
\displaystyle\frac{\sqrt{2a}}{\alpha}&0\\[3mm]
0&\displaystyle-\frac{\sqrt{2a}}{\alpha v_0-a}
\end{pmatrix}.
\]
Since $(\bar X_t)$ converges in distribution, it is straightforward to see that
$(h_t)$ converges in distribution to $0$. Concerning the Brownian martingale \(M\) its 
quadratic covariation matrix is given by
\[
\langle M\rangle_t=
\begin{pmatrix}
\displaystyle\int_0^t\bar X_s\,ds&t\\[3mm]
t&\displaystyle\int_0^t\frac{1}{\bar X_s}\,ds
\end{pmatrix}.
\]
Recalling the definitions of \(\bar I_t\) and \(\bar J_t\), we finally obtain
\[
\frac{1}{t}\langle M\rangle_t
=
\begin{pmatrix}
\bar I_t&1\\
1&\bar J_t
\end{pmatrix}
\xrightarrow[t\to\infty]{}
\Sigma=
\begin{pmatrix}
\displaystyle v_0&1\\[2mm]
\displaystyle 1&\frac{\alpha}{\alpha v_0-a}
\end{pmatrix}.
\]
Moreover, the matrix $\Sigma$ is positive definite since
$$
\det(\Sigma)=\frac{a}{\alpha v_0-a}>0.
$$
Since $(M_t)_{t\geq 0}$ is a continuous local martingale and
\[
\frac{1}{t}\langle M\rangle_t
\xrightarrow[t\to\infty]{\mathbb P}\Sigma,
\]
where $\Sigma$ is a deterministic positive definite matrix, the martingale central limit theorem (See Theorem \ref{Zanten2000}) yields
\[
\frac{M_t}{\sqrt{t}}
\xrightarrow[t\to\infty]{\mathcal L}
\mathcal N_2(0,\Sigma).
\]
In particular, the limiting Gaussian vector is non-degenerate.
By the continuous mapping theorem,
we obtain
\[
Q\frac{M_t}{\sqrt{t}}
\xrightarrow[t\to\infty]{\mathcal L}
\mathcal N_2(0,Q\Sigma Q^\top).
\]
A straightforward computation yields
\[
Q\Sigma Q^\top
=
\begin{pmatrix}
\displaystyle\frac{2av_0}{\alpha^2}
&
\displaystyle-\frac{2a}{\alpha(\alpha v_0-a)}
\\[3mm]
\displaystyle-\frac{2a}{\alpha(\alpha v_0-a)}
&
\displaystyle\frac{2a\alpha}{(\alpha v_0-a)^3}
\end{pmatrix}.
\]
Therefore,
\begin{equation}
\label{erreur_normale}
\sqrt{t}
\begin{pmatrix}
\bar I_t-v_0\\[2mm]
\displaystyle\bar J_t-\frac{\alpha}{\alpha v_0-a}
\end{pmatrix}
\xrightarrow[t\to\infty]{\mathcal L}
\mathcal N_2\left(0,Q\Sigma Q^\top\right).
\end{equation}
Now, let $h:\mathbb{R}^2\rightarrow\r$ be $C^3_b$, and, for $T>0$,
$$g_T:(x,y,z)\in\mathbb{R}^3\longmapsto h\left(\sqrt{T}\left(y - v_0\right),\sqrt{T}\left(z - \frac{1}{v_0-a/\alpha}\right)\right),$$
where~$x$ is a dummy variable. Then
\begin{align}\
\esp{g_T(0,I^N_T,J^N_T)} &= \left(\esp{g_T(0,I^N_T,J^N_T)}-\esp{g_T(0,\bar I_T,\bar J_T)}\right) + \esp{g_T(0,\bar I_T,\bar J_T)}\nonumber\\
&=\left(\esp{g_T(X^N_T,I^N_T,J^N_T)}-\esp{g_T(\bar X_T,\bar I_T,\bar J_T)}\right) + \esp{g_T(0,\bar I_T,\bar J_T)},\label{eq:ijan}
\end{align}
where the second term in~\eqref{eq:ijan} converges to
$\mathbb{E}[h(G)]$, with $G \sim \mathcal{N}_2(0, Q\Sigma Q^\top)$,
in virtue of relation~\eqref{erreur_normale}. To handle the first term of the RHS of~\eqref{eq:ijan}, notice that each function~$g_T$ is of class~$C^3_b$ and satisfies
$$\left|\left|g_T\right|\right|_{3,\infty} \leq C\, T^{3/2}\, \left|\left|h\right|\right|_{3,\infty},$$
with $C>0$ independent from~$T,h$. Then, re-using~\eqref{eq:thm21} with $g_T$ instead of~$g$, we obtain
$$\left|\mathbb{E}\left[g_T\left(X_T^N,I_T^N,J_T^N\right)\right] - \mathbb{E}\left[g_T\left(\bar X_T,\bar I_T,\bar J_T\right)\right]\right| \leq C ||h||_{3,\infty} T^{3/2}\left(N^{-1/2} e^{C_2\,T/\eps} + T\,\eps + \frac{\gamma}{T} + \eta_N\right).$$

Then, as in the end of the proof of Theorem~\ref{T1}, we choose successively the dependencies between the parameters~$T,N,\eps,\gamma,\eta_N$. The critical choice does not depend on~$\gamma,\eta_N$. For~$\eps,$ it is require that
$$T^{5/2}\eps_T$$
vanishes as $T$ goes to infinity, and~$N$ has to be chosen such that
$$N_T^{-1/2}\, T^{3/2} e^{C_2\,T/\eps_T}$$
goes to zero. This is always possible assuming that $N_T\geq e^{T^{7/2}M_T}$ where $M_T$ goes to infinity as $T$ goes to infinity.

This proves the convergence of $(I^N_T,J^N_T)$ for test-functions of class $C^3_b$. Thanks to Lemma~\ref{lem:cvLoi}, this implies the convergence in distribution.
\end{proof}

Let us now prove the asymptotic normality for the triple $(\widehat{v}^N_0(T),\widehat{a}^N(T),\widehat{\alpha}^N(T))$.
\begin{proof}[Proof of Theorem~\ref{T3}]
Firstly, from~\eqref{eq:a0} and~\eqref{a1}, it is clear that
$$\sqrt{T}\left(\widehat a^N(T) - a\right)\longrightarrow 0$$
in probability as $N,T$ tend jointly to infinity. Then, recalling that
\begin{multline*}
\left(\sqrt{T}\left(\widehat{v}^N_0(T) - v_0\right),\sqrt{T}\left(\widehat{a}^N(T) - a\right),\sqrt{T}\left(\widehat{\alpha}^N(T) - \alpha\right)\right) \\
= \left(\sqrt{T}\left(I^N_T - v_0\right),\sqrt{T}\left(\widehat{a}^N(T) - a\right),\sqrt{T}\left(\frac{\widehat{a}^N(T)}{I^N_T - 1/J^N_T} - \alpha\right)\right)
\end{multline*}
the statement of Theorem~\ref{T3} is straightforward by Slutsky's theorem and Proposition~\ref{prop:IJan}. Indeed, it is clear for the two first coordinates, and for the third one, it comes from the following calculation: denoting $\beta= v_0 - a/\alpha,$
\begin{align*}
    \frac{\widehat{a}^N(T)}{I^N_T - 1/J^N_T} - \alpha=& \frac{\widehat{a}^N(T)}{I^N_T - 1/J^N_T} - \frac{a}{v_0-\beta}\\
    =& \frac{\widehat{a}^N(T)(v_0-\beta) - a(I^N_T - 1/J^N_T)}{(I^N_T - 1/J^N_T)(v_0-\beta)}\\
    =& \frac{1}{(I^N_T - 1/J^N_T)(v_0-\beta)}\left[(\widehat{a}^N(T) - a)(v_0-\beta) + a(v_0 - I^N_T + (1/J^N_T)-\beta)\right]\\
    =&\frac{1}{(I^N_T - 1/J^N_T)(v_0-\beta)}\left[(\widehat{a}^N(T) - a)(v_0-\beta) - a(I^N_T - v_0) - \frac{a}{J^N_T/\beta}\left(J^N_T - \frac{1}{\beta}\right)\right],
\end{align*}
where Proposition~\ref{prop:IJan}, Theorem~\ref{T2} and Slutsky's theorem allow to conclude, recalling that $I^N_t$ (resp. $J^N_t$) converges in probability to $v_0$ (resp. $1/\beta$).
\end{proof}

\appendix
\section{Semigroups and generators}\label{sec:A}

The aim of this section is to state some technical lemmas about semigroups, that are used in our main results. All this section strongly relies on Appendices~A,~B and~E of \cite{erny2022mean}.

\begin{definition}\label{def:semi}
Let $(X_t)_{t\geq 0}$ be a $\r^d$-valued non-homogeneous Markov process (with $d\in\mathbb{N}_+$), and $(P_{s,t})_{s\leq t}$ be its semigroup: for any measurable function~$g:\r^d\to\r$ for which it makes sense,
$$P_{s,t} g(x) = \mathbb{E}_{\{X_s=x\}}\left[g(X_t)\right].$$
A measurable function~$g:\r^d\to\r$ is said to belong to the domain of the generator of~$X$ if there exists $(A_sg)_{s\geq 0}$ such that $t\mapsto P_{s,t}A_sg(x)$ is continuous at~$t=s$, and for all $x\in\r^d,s\leq t,$
\begin{itemize}
    \item $\displaystyle\mathbb{E}_{\{X_s=x\}}\left[g(X_t)\right] - g(x)= \mathbb{E}_{\{X_s=x\}} \int_s^t A_rg(X_r)dr,$
    \item $\displaystyle\mathbb{E}_{\{X_s=x\}} \int_s^t \left|A_rg(X_r)\right|dr<\infty.$
\end{itemize}
\end{definition}

Notice that the previous definition implies the classical property
$$\frac{d}{dt}\left(P_{s,t}g(x)\right)_{t=s} =A_sg(x).$$
And more generally, we have, respectively, the forward and backward equations
\begin{equation}\label{eq:forback}
\frac{d}{dt}\left(P_{s,t}g(x)\right) =P_{s,t}A_tg(x)\textrm{ and }\frac{d}{ds}\left(P_{s,t}g(x)\right) =-A_sP_{s,t}g(x).
\end{equation}
These equations can be proved with same techniques as in the proof of Proposition~A.3 of \cite{erny2022mean}.

A practical way to obtain an explicit form for a generator consists in using Ito's formula to write:
$$g(X_t)=g(X_s) + \int_s^t h_g(r,X_r)dr + M_{s,t},$$
with $(M_{s,t})_{t\geq s}$ a martingale null at $t=s$. This writing implies that $A_sg(x)=h_g(s,x)$.

The first technical result needed is the following.
\begin{lemma}\label{lem:trotter}
Let $(P^N_{r,t})_{r\leq t}$ and $(\bar P_{r,t})_{r\leq t}$
denote the (non-homogeneous) Markov semigroups associated respectively to $(Y^{N,\eps},I^{N,\gamma,\eps},J^{N,\gamma,\eps})$ and $(\bar X^\eps,\bar I^{\gamma,\eps},\bar J^{\gamma,\eps})$ defined in Section~\ref{sec:31}, and $(A^N_r)_{r\geq 0},(\bar A_r)_{r\geq 0}$ be their generator. Then for any~$g$ in $C^3_b(\r^3)$,
$$P^N_{s,t}g(x) - \bar P_{s,t}g(x) = \int_s^t P^N_{s,r}\left(A^N_r - \bar A_r\right)\bar P_{r,t}g(x)dr.$$
\end{lemma}

\begin{proof}[Sketch of proof]
Introduce
$$u:r\in[s,t]\longmapsto P^N_{s,r}\bar P_{r,t}g(x).$$

According to~\eqref{eq:forback}, $u$ is differentiable and, for $r\in(s,t)$
\begin{align*}
    u'(r) =& \frac{d}{dv}\left(P^N_{s,v}\bar P_{r,t}\right)_{v=r} + \frac{d}{dv}\left(P^N_{s,r}\bar P_{v,t}\right)_{v=r}\\
    =& P^N_{s,r}A^N_r\bar P_{r,t}g(x) - P^N_{s,r}\bar A_r\bar P_{r,t}g(x).
\end{align*}

The most technical part of the proof consists in showing that $u'$ is continuous, which can be achieved following the proof of Proposition~B.2 of \cite{erny2022mean}. Admitting this property, the fundamental theorem of calculus implies
$$u(t) - u(s) = \int_s^t u'(r)dr,$$
which is exactly the formula of the lemma.
\end{proof}

The second technical result that is required is the differentiability of the semigroup~$\bar P$ introduced above.

\begin{lemma}\label{lem:diffsemi}
Let $(\bar P_{r,t})_{r\leq t}$ be as in Lemma~\ref{lem:trotter}. For all $g\in C^3_b(\r^3)$, for any $y,z\in\r$ and $s\leq t$, the function
$$x\in\r\longmapsto \bar P_{s,t}g(x,y,z)$$
is $C^3_b$ with
$$\underset{x>0}{\sup}\left|\partial^3_x \bar P_{s,t} g(x,y,z)\right|\leq C_1 ||g||_{3,\infty} e^{C_2 t/\eps}$$
independently of $y,z.$
\end{lemma}

\begin{proof}[Sketch of proof]
In order to simplify the notation, we just prove the result for~$s=0$, and for the first-order derivative. Let us introduce the stochastic flow related to the process:
\begin{align*}
    \bar X^{(x)}_t =& x + \int_0^t b(\bar X^{(x)}_s)ds + \int_0^t \sigma(\bar X^{(x)}_s)dB_s,\\
    \bar I_t^{(x,y)} =& y + \frac{1}{t+\gamma}\int_0^t \bar X^{(x)}_sds,\\
    \bar J_t^{(x,z)} =& z + \frac{1}{t+\gamma}\int_0^t \frac{1}{\delta_\eps(\bar X^{(x)}_s) + \eta_N}ds.
\end{align*}
with $b(x) = \alpha(v_0-x)$ and $\sigma(x)=\sqrt{2a\delta_\eps(x)}$. The notation $\gamma,\eps$ in superscript are omitted in this proof.

Formally,
\begin{equation}\label{eq:semiflow}
\bar P_{0,t}g(x,y,z) = \mathbb{E}\left[g\left(\bar X^{(x)}_t,\bar I_t^{(x,y)},\bar J_t^{(x,z)}\right)\right].
\end{equation}
Hence, in order to differentiate the above expression w.r.t.~$x$, it is sufficient to differentiate (almost surely) the stochastic flow and to prove estimates for it.

Firstly, by Theorem~4.6.5 of \cite{kunita}, almost surely, for all~$t\geq 0$, the function $x\mapsto \bar X^{(x)}_t$ is differentiable, and so are the functions $x\mapsto \bar I^{(x,y)}_t$ and $x\mapsto \bar J^{(t,x)}_t$. It can be shown that
\begin{align*}
\partial_x \bar X^{(x)}_t =& 1 + \int_0^t \partial_x X^{(x)}_s b(X^{(x)}_s)ds + \int_0^t \partial_x X^{(x)}_s \sigma(X^{(x)}_s)dB_s.
\end{align*}

Then, using our estimates for $\bar X$ from Lemma~\ref{L1}, Burkholder-Davis-Gundy's inequality and Gr\"onwall's inequality, entails that: for any~$p\geq 1$ and~$t,x>0$,
\begin{equation}\label{eq:estimated1}
\mathbb{E}\left[\underset{s\leq t}{\sup}~|\partial_x \bar X^{(x)}_s|^p\right]\leq C_{p,t}\left(1+x^p\right).
\end{equation}

Whence, by the dominated convergence theorem (applied almost surely),
\begin{align*}
    \partial_x \bar I^{(x,y)}_t=& \frac{1}{t+\gamma}\int_0^t \partial_x\bar X^{(x)}_sds,\\
    \partial_x\bar J_t^{(x,z)} =& -\frac{1}{t+\gamma}\int_0^t \delta_\eps'(\bar X^{(x)}_s)\partial_x\bar X^{(x)}_s\frac{1}{\left(\delta_\eps(\bar X^{(x)}_s) + \eta_N\right)^2}ds,
\end{align*}
and by Jensen's inequality and H\"older's inequality, for any~$\mu>\lambda>1$,
\begin{align*}
    \mathbb{E}\left[\left|\partial_x \bar I^{(x,y)}_t\right|^\lambda\right]\leq & \frac{1}{t+\gamma}\int_0^t \mathbb{E}\left[\left|\partial_x\bar X^{(x)}_s\right|^\lambda\right]ds,\\
    \mathbb{E}\left[\left|\partial_x\bar J_t^{(x,z)}\right|^\lambda\right] \leq& \frac{1}{t+\gamma}\int_0^t \mathbb{E}\left[\left|\delta_\eps'(\bar X^{(x)}_s)\partial_x\bar X^{(x)}_s\frac{1}{\left(\delta_\eps(\bar X^{(x)}_s) + \eta_N\right)^2}ds\right|^\lambda\right]\\
    \leq& \frac{1}{t+\gamma}\int_0^t \mathbb{E}\left[\left|\partial_x\bar X^{(x)}_s\right|^{\mu/(\mu-\lambda)}\right]^{(\mu-\lambda)/\mu}\mathbb{E}\left[\left|{\left(\bar X^{(x)}_s + \eta_N\right)}\right|^{-2\mu}\right]^{\lambda/\mu}ds\\
    \leq& \left(\eps/2\right)^{-2\lambda}\frac{1}{t+\gamma}\int_0^t \mathbb{E}\left[\left|\partial_x\bar X^{(x)}_s\right|^{\mu/(\mu-\lambda)}\right]^{(\mu-\lambda)/\mu}ds,
\end{align*}
where we have used that $|\delta_\eps'|$ is bounded by one, and that $\delta_\eps$ is lower-bounded by $\eps/2$.

The above inequalities together with~\eqref{eq:estimated1} imply that the semigroup~\eqref{eq:semiflow} is differentiable w.r.t.~$x$ with
\begin{align*}
    \partial_x \bar P_{0,t}g(x,y,z) =& \mathbb{E}\left[(\partial_x \bar X^{(x)}_t) \partial_1 g\left(\bar X^{(x)}_t,\bar I_t^{(x,y)},\bar J_t^{(x,z)}\right)\right]\\
    &+\mathbb{E}\left[(\partial_x \bar I^{(x,y)}_t) \partial_2 g\left(\bar X^{(x)}_t,\bar I_t^{(x,y)},\bar J_t^{(x,z)}\right)\right]\\
    &+\mathbb{E}\left[(\partial_x \bar J^{(x,z)}_t) \partial_3 g\left(\bar X^{(x)}_t,\bar I_t^{(x,y)},\bar J_t^{(x,z)}\right)\right].
\end{align*}

Then, classical argument of stochastic analysis and Gr\"onwall's inequality leads to
$$\underset{x>0}{\sup}\left|\partial^1_x \bar P_{0,t} g(x,y,z)\right|\leq C_3(t,\eps) ||g||_{1,\infty} e^{C_4(\eps) t}.$$

We refer to the calculations in the proof of Proposition~3.4 (in Appendix~E) of \cite{erny2022mean} for the form of the above functions~$C_3,C_4$: there exist positive constants~$C_1,C_2$ and $n\in\mathbb{N}^*,$
$$C_3(t,\eps)\leq C_1 \left((1+(t/\eps)^n\right)\textrm{ and }C_4(\eps) \leq C_2 /\eps.$$
Notice that the dependencies of $C_4(\eps)$ w.r.t.~$1/\eps$ above comes from the fact that the Lipschitz constant of the coefficient~$\sigma(x)=\sqrt{2a \delta_\eps(x)}$ is non-greater than $C/\sqrt{\eps}$.

Finally
$$\underset{x>0}{\sup}\left|\partial^1_x \bar P_{s,t} g(x,y,z)\right|\leq C_1\left(1+(t/\eps)^n\right) ||g||_{1,\infty} e^{C_2 t/\eps}\leq C_1 e^{C_2 t/\eps},$$
where the last inequality is for free by possibly increasing the constant~$C_2$. This proves the desired control for the first-order derivative. For the second-order and the third-order ones, similar calculations lead to the same kind of control (the values of $C_1,C_2,C_3,C_4,n$ above change, but the overall dependencies w.r.t.~$t,\eps$ remain the same). We refer again to the proof of Proposition~3.4 of \cite{erny2022mean} for these calculations in a simpler but similar model.
\end{proof}

\section{Proof of Lemma~\ref{L1}: a priori estimates}\label{sec:B}

For Item~$(i)$, the estimate for $\bar X$ is classical, and holds true more generally for any process solution to an SDE with sublinear coefficients. The estimate for $X^N$ is finer and require to exploit the compensation of the two jump terms for the control to be uniform w.r.t.~$N$ (similarly as in the proof of Lemma~3.1.(i) of \cite{erny2022mean} for example).

Now let us prove~$(ii)$. By definition
$$\bar X_t = x + \int_0^t \alpha(v_0 - \bar X_s)ds + \int_0^t \sqrt{2a \bar X_s}dB_s.$$
And, by $(i)$, the Brownian integral above is a (real) martingale, hence its expectation is zero. This implies that $t\mapsto \mathbb{E}_x[\bar X_t]$ solves the ODE $y' = \alpha(v_0-y)$ with $y(0)=x$, which entails the result. The same calculations hold for $X^N$.

Item~$(iii)$ is exactly Lemma~3.1.(i) of \cite{ben2020local}.

To prove~$(iv)$, let us denote, for $x>-1/\sqrt{N}$,
$$\phi_N^p(x) = (x + \eta_N)^{-p}.$$

Notice that
$$(\phi_N^p)' = -p\, \phi_N^{p+1}.$$

Then,
\begin{align*}
A^N\phi_N^p(x)=& \alpha(v_0-x)(\phi_N^p)'(x) + N\,a\, x_+\left[\phi_N^p(x+1/\sqrt{N}) - \phi_N^p(x) + \phi_N^p(x-1/\sqrt{N}) - \phi_N^p(x)\right]\\
=& p\alpha x\,\phi_N^{p+1}(x) - p\alpha v_0\, \phi_N^{p+1}(x) \\
&+ N\,a\, x_+\left[\phi_N^p(x+1/\sqrt{N}) - \phi_N^p(x) + \phi_N^p(x-1/\sqrt{N}) - \phi_N^p(x)\right]
\end{align*}

And, since
$$x\,\phi_N^{p+1(x)} = \frac{x}{(x+\eta_N)^{p+1}}\leq \frac{x+\eta_N}{(x+\eta_N)^{p+1}} = \phi_N^p(x),$$
we have
\begin{align*}
A^N\phi_N^p(x) \leq& p\alpha\,\phi_N^{p}(x) - p\alpha v_0\, \phi_N^{p+1}(x) \\
&+ N\,a\, x_+\left[\phi_N^p(x+1/\sqrt{N}) - \phi_N^p(x) + \phi_N^p(x-1/\sqrt{N}) - \phi_N^p(x)\right]
\end{align*}

Furthermore, by Taylor-Lagrange's formula,
\begin{align*}
&\phi_N^p(x+1/\sqrt{N}) - \phi_N^p(x) + \phi_N^p(x-1/\sqrt{N}) - \phi_N^p(x)\\
&=\frac1N(\phi_N^p)''(x) + \frac{1}{6N^{3/2}}\left[(\phi_N^p)'''(u) - (\phi_N^p)'''(v)\right]\\
&=\frac1N(\phi_N^p)''(x) + \frac{1}{6N^{3/2}}p(p+1)(p+2)\left[(\phi_N^{p+3})(v) - (\phi_N^{p+3})(u)\right]
\end{align*}
for some
$$x-\frac1{\sqrt{N}}\leq v\leq x\leq u \leq x+\frac{1}{\sqrt{N}}.$$

The function $\phi_N^{p+3}$ being non-increasing, the expression in the brackets is positive, and, by positivity of the function,
\begin{align}
&\phi_N^p(x+1/\sqrt{N}) - \phi_N^p(x) + \phi_N^p(x-1/\sqrt{N}) - \phi_N^p(x)\nonumber\\
&\leq \frac1N(\phi_N^p)''(x) + \frac{1}{6N^{3/2}}p(p+1)(p+2)\phi_N^{p+3}(x-1/\sqrt{N})\label{eq:contSaut}
\end{align}

Now, fix some $\lambda >0$ (whose value will be chosen arbitrarily close to zero later), and notice that, for $N$ large enough and for any~$x>-1/\sqrt{N}$
$$\frac{1}{\sqrt{N}}\leq \lambda(\eta_N - 2/\sqrt{N})\leq \lambda (x + \eta_N - 1/\sqrt{N}).$$

This implies that
$$\frac{1}{\sqrt{N}}\phi_N^{p+3}(x-1/\sqrt{N})\leq \lambda \phi_N^{p+2}(x-1/\sqrt{N})\leq 2\lambda \phi_N^{p+2}(x).$$

Using the inequality above in~\eqref{eq:contSaut} leads to
\begin{align}
&\phi_N^p(x+1/\sqrt{N}) - \phi_N^p(x) + \phi_N^p(x-1/\sqrt{N}) - \phi_N^p(x)\nonumber\\
&\leq \frac1N p(p+1)\phi_N^{p+2}(x) + \frac{1}{3 N}p(p+1)(p+2)\lambda\phi_N^{p+2}(x)\nonumber\\
&\leq \frac1N p(p+1)\phi_N^{p+2}(x)\left[1 + \frac{\lambda}{3}(p+2)\right].\nonumber
\end{align}

Coming back to the generator,
\begin{align*}
A^N \phi_N^p(x)\leq& p\alpha\,\phi_N^p(x) - p\alpha v_0\,\phi_N^{p+1}(x) + a\,x_+ p(p+1)\phi_N^{p+2}(x)\left[1 + \frac{\lambda}{3}(p+2)\right]
\end{align*}

Using that
$$x_+ \phi_N^{p+2}(x) \leq (x+\eta_N)\phi_N^{p+2}(x) = \phi_N^{p+1}(x),$$
we obtain the following control
$$A^N\phi_N^p(x) \leq p\alpha\,\phi_N^p(x) + p \phi_N^{p+1}(x)\left[a(p+1)\left(1+\frac{\lambda}{3}(p+2)\right) - \alpha v_0\right].$$

Then, recalling that $a(p+1)<\alpha v_0$, it is possible to chose a value for $\lambda>0$ such that the expression in the brackets above is a negative constant, implying
$$A^N\phi_N^p(x) \leq p\alpha\, \phi_N^p(x) - c_p \phi_N^{p+1}(x),$$
with $c_p>0.$

Then, denoting $y = x + \eta_N$, and fixing $\delta>0$ (that will be chosen arbitrarily small later),
$$\frac{1}{y^p} \leq \frac{\delta}{y^{p+1}} + \frac{1}{\delta^p}.$$

Indeed, if $y\geq \delta$ then $y^{-p}\leq \delta^p$, and otherwise $y^{-p} = y y^{-p-1}\leq \delta y^{-p-1}.$

Then
$$A^N \phi_N^p(x) \leq \frac{p\alpha}{\delta^p} + \phi_N^{p+1}(x)\left[\delta - c_p\right].$$

In particular, choosing $\delta$ smaller than $c_p$ implies that
\begin{equation}\label{eq:ANboudnd}
A^N\phi_N^p(x)\leq \frac{p\alpha}{\delta^p}.
\end{equation}

Finally, by Ito's formula,
$$\phi^p_N\left(X^N_t\right) = \phi_N^p\left(X^N_0\right) + \int_0^t A^N\phi_N^p\left(X^N_s\right)ds + M^N_t,$$
with $M^N_t$ the local martingale defined as
\begin{align*}
M^N_t =& \sum_{j=1}^N \int_0^t\int_0^\infty \uno{z\leq a (X^N_{s-})_+}\left[\phi_N^p(X^N_{s-} + 1/\sqrt{N}) - \phi_N^p(X^N_{s-})\right]d\tilde\pi_j^+(s,z) \\
&+ \sum_{j=1}^N \int_0^t\int_0^\infty \uno{z\leq a (X^N_{s-})_+}\left[\phi_N^p(X^N_{s-} - 1/\sqrt{N}) - \phi_N^p(X^N_{s-})\right]d\tilde\pi_j^-(s,z)
\end{align*}
where $\tilde \pi_j^+,\tilde \pi_j^-$ are the compensated versions of $\pi_j^+,\pi_j^-$ (for $1\leq j\leq N$). Since the function $\phi_N^p$ is bounded (by $(\eta_N-1/\sqrt{N})^p$), it is classical that $M^N$ is a true martingale, whence
\begin{align*}
\mathbb{E}_x\left[\phi^p_N\left(X^N_t\right)\right]=\phi_N^p(x) + \int_0^t \mathbb{E}_x\left[A^N\phi^p_N\left(X^N_s\right)\right]ds\leq \phi_N^p(x) + \frac{p\alpha}{\delta^p}t,
\end{align*}
which ends the proof.
\section{A Central Limit Theorem for Continuous Local Martingales}\label{sec:C}

The next theorem is about the asymptotic behaviour of continuous multivariate local martingales, see \cite[Theorem 4.1]{Zanten2000}.
\begin{theorem}\label{Zanten2000}
	For $p\in\mathbb{N} \setminus \lbrace 0\rbrace$, let $M=(M_t)_{t\in\R_+}$ be a $p$-dimensional square-integrable continuous
	local martingale with respect to the filtration $(\mathcal{F}_t)_{t\in\R_+}$ such that $\mathbb{P}(M_0=0)=1$. Suppose that
	there exists a function $Q:\left[t_0,\infty\right) \to \mathcal{M}_p$ with some $t_0\in\R_+$ such that $Q(t)$ is an invertible
	(non-random) matrix for all $t\in\R_+$, $\underset{t\to\infty}{\lim} ||Q(t)||=0$ and
	$$	Q(t) \langle M\rangle_t Q(t)^{\mathsf{T}}\stackrel{\mathbb{P}}{\longrightarrow}\eta\eta^{\mathsf{T}},\quad \text{as }t\to \infty,$$
	where $\eta$ is a random matrix in $\mathcal{M}_p$. Then, for each random matrix $A\in\mathcal{M}_{k,l}$,  $k,l\in\mathbb{N}\setminus\lbrace 0\rbrace$, defined on $(\Omega,\mathcal{F},\mathcal{P})$, we have
	$$
	(Q(t)M_t,A) \stackrel{\mathcal{D}}{\longrightarrow}
	(\eta Z,A),\quad \text{as }t\to \infty,$$
	where $Z$ is a $p$-dimensional standard normally distributed random vector independent of $(\eta,A)$.
	\label{CLT Van Zanten}
\end{theorem}
\bibliography{ref} 
\end{document}